\input ./style/arxiv-general.cfg
\documentclass[aop,MSNbibl,dvips]{arximspdf}
\makeatletter
   \@ifpackageloaded{graphicx}{}{\usepackage{graphicx}}
\makeatother

\doi{10.1214/15-AOP1025}
\volume{44}
\issue{4}
\pubyear{2016}
\firstpage{2479}
\lastpage{2506}
\docsubty{FLA}

\makeatletter
\newcommand{\lleft}{\left}
\newcommand{\rrvert}{\vert}
\newcommand{\rright}{\right}
\newcommand{\rrVert}{\Vert}
\newcommand{\llvert}{\vert}
\newcommand{\llVert}{\Vert}
\renewcommand{\mid}{|}
\newtheorem{teo}{Theorem}[section]
\newtheorem{lem}[teo]{Lemma}
\newtheorem{prop}[teo]{Proposition}
\newtheorem{cor}[teo]{Corollary}
\newproclaim{rem}[teo]{Remark}
\newcommand{\Var}{\operatorname{Var}}
\newcommand{\EE}{\mathbf{E}}
\newcommand{\Tr}{\operatorname{Tr}}
\makeatother

\begin{document}
\begin{frontmatter}

\title{Sharp nonasymptotic bounds on the norm of random matrices with
independent entries}
\runtitle{The norm of random matrices with independent entries}

\begin{aug}
\author[A]{\fnms{Afonso S.}~\snm{Bandeira}\ead[label=e1]{ajsb@math.princeton.edu}\thanksref{T1}}
\and
\author[B]{\fnms{Ramon}~\snm{van Handel}\corref{}\ead[label=e2]{rvan@princeton.edu}\thanksref{T2}}
\runauthor{A.~S. Bandeira and R. van Handel}
\affiliation{Princeton University}
\address[A]{Program in Applied\\
\quad and Computational Mathematics\\
Princeton University\\
Princeton, New Jersey 08544\\
USA\\
\printead{e1}}
\address[B]{Department of Operations Research\\
\quad and Financial Engineering\\
Sherrerd Hall 227\\
Princeton University \\
Princeton, New Jersey 08544\\
USA\\
\printead{e2}}
\end{aug}
\thankstext{T1}{Supported by AFOSR Grant FA9550-12-1-0317.}
\thankstext{T2}{Supported in part by NSF Grant CAREER-DMS-1148711 and
by the ARO through PECASE award W911NF-14-1-0094.}

%
\received{\smonth{8} \syear{2014}}
%
\revised{\smonth{3} \syear{2015}}

%
\begin{abstract}
We obtain nonasymptotic bounds on the spectral norm of random
matrices with independent entries that improve significantly on earlier
results. If $X$ is the $n\times n$ symmetric
matrix with $X_{ij}\sim N(0,b_{ij}^2)$, we show that
\[
\mathbf{E}\llVert X\rrVert\lesssim\max_i\sqrt{
\sum_{j}b_{ij}^2} + \max
_{ij}\llvert b_{ij}\rrvert\sqrt{\log n}.
\]
This bound is optimal in the sense that a matching lower bound holds under
mild assumptions, and the constants are sufficiently sharp that we can
often capture the precise edge of the spectrum. Analogous
results are obtained for rectangular matrices and for more general
sub-Gaussian or heavy-tailed distributions of the entries, and we derive
tail bounds in addition to bounds on the expected norm. The proofs are
based on a combination of the moment method and geometric functional
analysis techniques.
As an application, we show that our bounds immediately yield the correct
phase transition behavior of the spectral edge of random band matrices and
of sparse Wigner matrices. We also recover a result of Seginer on
the norm of Rademacher matrices.
\end{abstract}

%
\end{frontmatter}

\section{Introduction}\label{sec1}

Understanding the behavior of the spectral norm of random matrices is a
fundamental problem in probability theory, as well as a problem of
considerable importance in many modern applications. If $X$ is a Wigner
matrix, that is, a symmetric $n\times n$ matrix whose entries are
i.i.d. with unit variance, then a classical result in random matrix theory
\cite{FK81,BY88,AGZ10,Tao12} shows that $\llVert X\rrVert /\sqrt{n}\to
2$ under mild
moment
assumptions (as is expected from the well-known fact that the empirical
spectral density converges to the semicircle law supported in $[-2,2]$).
The corresponding result for rectangular matrices with i.i.d. entries is
even older~\cite{Gem80}. More recently, there has been considerable
interest in structured random matrices where the entries are no longer
identically distributed. As the combinatorial methods that are used for
this purpose typically exploit the specific structure of the entries,
precise asymptotic results on the spectral norm of structured matrices
must generally be obtained on a case-by-case basis; see, for example,
\cite{Sod09,Sod10}.

In order to gain a deeper understanding of the spectral norm of structured
matrices, it is natural to ask whether one can find a unifying principle
that captures at least the correct scale of the norm in a general setting,
that is, in the absence of specific structural assumptions. This question
is most naturally phrased in a nonasymptotic setting: can we obtain upper
and lower bounds on $\llVert X\rrVert $, in terms of natural
parameters that capture
the structure of $X$, that differ only by universal constants?
Nonasymptotic bounds on the norm of a random matrix have long been
developed in a different area of probability that arises from problems in
geometric functional analysis, and have had a significant impact on
various areas of pure and applied mathematics
\cite{DS01,RV10,Ver12}. Unfortunately, as we will shortly see, the
best known general results along these lines fail to capture the correct
scale of the spectral norm of structured matrices except in extreme
cases.

In this paper, we investigate the norm of random matrices with independent
entries. Consider for concreteness the case of Gaussian matrices
(our main results will extend to more general distributions of the
entries). Let $X$ be the $n\times n$ symmetric random matrix with entries
$X_{ij}=g_{ij}b_{ij}$, where $\{g_{ij}\dvtx i\ge j\}$ are independent standard
Gaussian random variables and $\{b_{ij}\dvtx i\ge j\}$ are given scalars.

Perhaps the most useful known nonasymptotic bound on the spectral norm
$\llVert X\rrVert $ can be obtained as a consequence of the
noncommutative Khintchine
inequality of Lust-Piquard and Pisier \cite{Pis03}, or alternatively
(in a much more elementary fashion) from the ``matrix concentration''
method that has been widely developed in recent years \cite{Oli10,Tro12}.
This yields the following inequality in our setting:
\[
\EE\llVert X\rrVert\lesssim\sigma\sqrt{\log n} \qquad\mbox{with }
\sigma:=
\max_i\sqrt{\sum_j
b_{ij}^2}. %
\]
Unfortunately, this inequality already fails to be sharp in the simplest
case of Wigner matrices: here $\sigma=\sqrt{n}$, so that the resulting
bound $\EE\llVert X\rrVert \lesssim\sqrt{n\log n}$ falls short of the
correct scaling
$\EE\llVert X\rrVert \sim\sqrt{n}$. On the other hand, the
logarithmic factor in this
bound is necessary: if $X$ is the diagonal matrix
with independent standard Gaussian entries, then $\sigma=1$ and
$\EE\llVert X\rrVert \sim\sqrt{\log n}$. We therefore
conclude that while the noncommutative Khintchine bound is sharp in
extreme cases, it fails to capture the structure of the matrix $X$ in a
satisfactory manner.

A different bound on $\llVert X\rrVert $ can be obtained by a method
due to Gordon
(see \cite{DS01}) that exploits Slepian's comparison lemma for Gaussian
processes, or alternatively from a simple $\varepsilon$-net argument
\cite{Ver12,Tao12}. This yields the following inequality:
\[
\EE\llVert X\rrVert\lesssim\sigma_*\sqrt{n}\qquad\mbox{with }
\sigma_*:=
\max_{ij}\llvert b_{ij}\rrvert. %
\]
While the parameter $\sigma_*$ that appears in this bound is often much
smaller than $\sigma$, the dimensional scaling of this bound is much worse
than in the noncommutative Khintchine bound. In particular, while this
bound captures the correct $\sqrt{n}$ rate for Wigner matrices, it is
vastly suboptimal in almost every other situation (e.g., in the
diagonal matrix example considered above).

Further nonasymptotic bounds on $\llVert X\rrVert $ have been
obtained in the present
setting by Lata{\l}a \cite{Lat05} and by Riemer and Sch\"utt \cite{RS13}.
In most examples, these bounds provide even worse rates than the
noncommutative Khintchine bound. Seginer \cite{Seg00} obtained a slight
improvement on the noncommutative Khintchine bound that is specific to the
special case where the random matrix has uniformly bounded entries; see
Section~\ref{secrad} below. None of these results provides a sharp
understanding of the scale of the spectral norm for general structured
matrices.

The present paper develops a new family of nonasymptotic bounds on the
spectral norm of structured random matrices that prove to be optimal in a
surprisingly general setting. Our main bounds are of the form
\[
\EE\llVert X\rrVert\lesssim\sigma+ \sigma_*\sqrt{\log n}, %
\]
which provides a sort of interpolation between the two bounds discussed
above. For example, the following is one of the main results of this
paper.

%
\begin{teo}
\label{teomain}
Let $X$ be the $n\times n$ symmetric matrix with
$X_{ij}=g_{ij}b_{ij}$, where $\{g_{ij}\dvtx i\ge j\}$ are i.i.d. $\sim N(0,1)$
and $\{b_{ij}\dvtx i\ge j\}$ are given scalars. Then
\[
\EE\llVert X\rrVert\le(1+\varepsilon)\biggl\{2\sigma+ \frac{6}{\sqrt
{\log(1+\varepsilon)}}
\sigma_*\sqrt{\log n}\biggr\} %
\]
for any $0<\varepsilon\le1/2$, where $\sigma,\sigma_*$ are as defined
above.
\end{teo}

\noindent
Let us emphasize two important features of this result:
\begin{itemize}
\item It is almost trivial to obtain a matching lower bound of the form
\[
\EE\llVert X\rrVert\gtrsim\sigma+\sigma_*\sqrt{\log n} %
\]
that holds as long as the coefficients $b_{ij}$ are not too inhomogeneous
(Section~\ref{seclower}). This means that Theorem~\ref{teomain}
captures the optimal scaling of the expected norm $\EE\llVert
X\rrVert $ under
surprisingly minimal structural assumptions.
\end{itemize}

\begin{itemize}
\item In the case of Wigner matrices, Theorem~\ref{teomain} yields a
bound of the form
\[
\EE\llVert X\rrVert\le(1+\varepsilon)2\sqrt{n}+o(\sqrt{n}) %
\]
for arbitrarily small $\varepsilon>0$. Thus Theorem~\ref{teomain} not
only captures the correct scaling of the spectral norm, but even recovers
the precise asymptotic behavior $\llVert X\rrVert /\sqrt{n}\to2$ as
$n\to\infty$.
This feature of Theorem~\ref{teomain} makes it possible to effortlessly
prove nontrivial results, such as the precise phase transition behavior of
the spectral edge of random band matrices (Section~\ref{secsparse}), that
would be distinctly nontrivial to obtain by classical combinatorial
methods.
\end{itemize}
In view of these observations, it seems that Theorem~\ref{teomain} is
essentially the optimal result of its kind: there is little hope to
accurately capture inhomogeneous models where Theorem~\ref{teomain} is
not sharp in terms of simple parameters such as $\sigma,\sigma_*$; see
Remarks~\ref{remdimfree} and \ref{remschutt}. On the other hand,
we can now understand the previous bounds as extreme cases of
Theorem~\ref{teomain}. The noncommutative Khintchine bound matches
Theorem~\ref{teomain} when $\sigma/\sigma_*\lesssim1$: this case is
minimal as $\sigma/\sigma_*\ge1$. Gordon's bound matches
Theorem~\ref{teomain} when $\sigma/\sigma_*\gtrsim\sqrt{n}$: this
case is maximal as $\sigma/\sigma_*\le\sqrt{n}$. In intermediate
regimes, Theorem~\ref{teomain} yields a strictly better scaling.

While we have formulated the specific result of Theorem~\ref{teomain} for
concreteness, our methods are not restricted to this particular setting.
Once a complete proof of Theorem~\ref{teomain} has been given in
Section~\ref{secmain}, we will develop various extensions and
complements in
Section~\ref{secexcom}. These results are developed both for their
independent interest, and in view of their potential utility in
applications to other problems. In Section~\ref{secrect}, we prove a
sharp analogue of Theorem~\ref{teomain} in the setting of rectangular
matrices. In Section~\ref{secnong}, we develop versions of our main
results when the entries are not necessarily Gaussian: for sub-Gaussian
variables we obtain very similar results to the Gaussian case, while the
scaling in our main results must be modified in the case of heavy-tailed
entries. In Section~\ref{secdimfree}, we develop variants of our main
results where the explicit-dimensional dependence is replaced by a certain
notion of effective dimension. In Section~\ref{sectail}, we obtain sharp
inequalities for the tail probabilities of the spectral norm $\llVert
X\rrVert $
rather than for its expectation. Finally, we obtain in
Section~\ref{seclower} lower bounds on the spectral norm of Gaussian matrices
that match our upper bounds under rather mild assumptions.

In Section~\ref{secex}, we develop two applications that illustrate the
power of our main results. In Section~\ref{secsparse}, we investigate a
phase transition phenomenon for the spectral edge of random band matrices
and more general sparse Wigner matrices. Our main results effortlessly
provide a precise understanding of this transition, which sharpens earlier
results that were obtained by much more delicate combinatorial methods in
\cite{Kho08,Sod10,BGP14}. In Section~\ref{secrad}, we investigate the
setting of Rademacher random matrices with entries
$X_{ij}=\varepsilon_{ij}b_{ij}$, where $\varepsilon_{ij}$ are independent
Rademacher (symmetric Bernoulli) variables. Here we recover a result of
Seginer \cite{Seg00} with a much simpler proof, and develop insight into
how such bounds can be improved.

One of the nice features of Theorem~\ref{teomain} is that its proof
explains very clearly why the result is true. Once the idea has been
understood, the technical details prove to be of minimal difficulty, which
suggests that the ``right'' approach has been found. Let us briefly
illustrate the idea behind the proof in the special case where the
coefficients $b_{ij}$ take only the values $\{0,1\}$ (this setting guided
our intuition, though the ultimate proof is no more difficult in the
general setting). We can then interpret the matrix of
coefficients $(b_{ij})$ as the adjacency matrix of a graph $G$ on $n$
points, and we have $\sigma_*=1$ and $\sigma=\sqrt{k}$ where $k$ is the
maximal degree of $G$.

Following\vspace*{1pt} a classical idea in random matrix theory, we use the fact that
the spectral norm $\llVert X\rrVert $ is comparable to the quantity
$\Tr[X^p]^{1/p}$
for $p\sim\log n$. If one writes out the expression for $\EE\Tr
[X^p]$ in
terms of the coefficients, it is readily seen that controlling this
quantity requires us to count the number of cycles in $G$ for which every
edge is visited an even number of times.
One might expect that the graph $G$ of degree $k$ that possesses the most
such cycles is the complete graph on $k$ points. If this were the case,
then one could control
$\EE\Tr[X^p]$ by $\EE\Tr[Y^p]$ where $Y$ is a Wigner matrix of dimension
$k$. This intuition is almost, but not entirely correct: while a
$k$-clique typically possesses more distinct topologies of cycles, each
cycle of a given topology can typically be embedded in more ways in a
regular graph on $n$ points than in a $k$-clique. Careful bookkeeping
shows that the latter can be accounted for by choosing a slightly larger
Wigner matrix of dimension $k+p$. We therefore obtain a comparison theorem
between the spectral norm of $X$ and the spectral norm of a
$(k+p)$-dimensional\vspace*{1pt} Wigner matrix, which is of the desired order
$\sqrt{k+p}\sim\sqrt{k}+\sqrt{\log n}$ for $p\sim\log n$. We can now
conclude by using standard ideas from probability in Banach spaces to
obtain sharp nonasymptotic bounds on the norm of the resulting Wigner
matrix, avoiding entirely any combinatorial complications. (A purely
combinatorial approach would be nontrivial as very high moments of Wigner
matrices can appear in this argument.)

We conclude the \hyperref[sec1]{Introduction} by noting that both the noncommutative
Khintchine inequality and Gordon's bound can be formulated in a more
general context beyond the case of independent entries. Whether the
conclusion of Theorem~\ref{teomain} extends to this situation is a
natural question of considerable interest.

\subsection*{Notation}

Let us clarify a few notational conventions that will be used throughout
the paper. In the sequel, $a\lesssim b$ means that $a\le Cb$ for a
universal constant $C$. (If $C$ depends on other quantities,
this will be indicated explicitly.) If $a\lesssim b$ and $b\lesssim a$, we
write $a\asymp b$. We write $a\wedge b:=\min(a,b)$ and $a\vee b
:=\max(a,b)$, and we denote by $[n]:=\{1,\ldots,n\}$. Finally, we
occasionally write $\llVert \xi\rrVert _p := \EE[\xi^p]^{1/p}$.

\section{Proof of Theorem~\texorpdfstring{\protect\ref{teomain}}{1.1}}\label{secmain}

The main idea behind the proof of Theorem~\ref{teomain} is the
following comparison theorem.

%
\begin{prop}
\label{propmaincomp}
Let $Y_r$ be the $r\times r$ symmetric matrix such that\break
$\{(Y_r)_{ij}\dvtx i\ge j\}$ are independent $N(0,1)$ random variables, and
suppose that $\sigma_*\le1$. Then
\[
\EE\Tr\bigl[X^{2p}\bigr] \le\frac{n}{\lceil\sigma^2\rceil+p} \EE\Tr
\bigl[Y^{2p}_{\lceil\sigma^2\rceil+p} \bigr] \qquad\mbox{for every }p\in
\mathbb{N}.
\]
\end{prop}

Let us begin by completing the proof of Theorem~\ref{teomain}
given this result. We need the following lemma, which is a variation
on standard ideas; cf. \cite{DS01}.

%
\begin{lem}
\label{lemgordon}
Let $Y_r$ be the $r\times r$ symmetric matrix such that
$\{(Y_r)_{ij}\dvtx i\ge j\}$ are independent $N(0,1)$ random variables. Then
for every $p\ge2$
\[
\EE\bigl[\llVert Y_r\rrVert^{2p}\bigr]^{1/2p} \le
2\sqrt{r} + 2\sqrt{2p}. %
\]
\end{lem}

\begin{pf}
We begin by noting that
\[
\llVert Y_r\rrVert=\lambda_+\vee\lambda_-,\qquad\lambda_+:=\sup
_{v\in S}\langle v,Y_rv\rangle, \qquad\lambda_-=-\inf
_{v\in S}\langle v,Y_rv\rangle, %
\]
where $S$ is the unit sphere in $\mathbb{R}^r$. We are therefore
interested in the supremum of the Gaussian process $\{\langle
v,Y_rv\rangle\}_{v\in S}$, whose natural distance can be estimated as
\[
\EE\bigl\llvert\langle v,Y_rv\rangle-\langle w,Y_rw
\rangle\bigr\rrvert^2 \le2\sum_{i,j}
\{v_iv_j-w_iw_j\}^2
\le4\llVert v-w\rrVert^2 %
\]
[using $1-x^2\le2(1-x)$ for $x\le1$]. The right-hand side of this
expression is the
natural distance of the Gaussian process $\{2\langle v,g\rangle\}
_{v\in
S}$, where $g$ is the standard Gaussian vector in $\mathbb{R}^r$.
Therefore, Slepian's lemma \cite{BLM13}, Theorem~13.3, implies
\[
\EE\lambda_+ = \EE\sup_{v\in S}\langle v,Y_rv\rangle
\le2\EE\sup_{v\in S}\langle v,g\rangle= 2\EE\llVert g\rrVert\le2
\sqrt{r}. %
\]
Moreover, note that $\lambda_+$ and $\lambda_-$ have the same
distribution (as evidently $Y_r$ and $-Y_r$ have the same distribution).
Therefore, using the triangle inequality for $\llVert \cdot\rrVert _{2p}$,
\begin{eqnarray*}
\EE\bigl[\llVert Y_r\rrVert^{2p}\bigr]^{1/2p} &=&
\llVert\lambda_+\vee\lambda_-\rrVert_{2p}
\\
&\le&\EE\lambda_+ + \llVert\lambda_+\vee\lambda_- - \EE\lambda_+\rrVert
_{2p}
\\
& =& \EE\lambda_+ + \bigl\llVert(\lambda_+-\EE\lambda_+)\vee(\lambda
_--\EE
\lambda_-)\bigr\rrVert_{2p}.
\end{eqnarray*}
It follows from
Gaussian concentration \cite{BLM13}, Theorems~5.8~and~2.1, that
\[
\EE\bigl[(\lambda_+-\EE\lambda_+)^{2p}\vee(\lambda_- -\EE
\lambda_-)^{2p}\bigr] \le p!4^{p+1} \le(2\sqrt{2p})^{2p}
\]
for $p\ge2$. Putting together the above estimates completes the proof.
\end{pf}

\begin{pf*}{Proof of Theorem~\ref{teomain}}
We can clearly assume without loss of generality that the matrix $X$ is
normalized such that $\sigma_*=1$. For $p\ge2$, we can estimate
\begin{eqnarray*}
\EE\llVert X\rrVert&\le&\EE\bigl[\Tr\bigl[X^{2p}\bigr]
\bigr]^{1/2p}
\\
&\le& n^{1/2p} \EE\bigl[\llVert Y_{\lceil\sigma^2\rceil+p}\rrVert^{2p}
\bigr]^{1/2p}
\\
&\le& n^{1/2p}\bigl\{2\sqrt{\bigl\lceil\sigma^2\bigr
\rceil+p}+2\sqrt{2p}\bigr\}
\end{eqnarray*}
by Proposition~\ref{propmaincomp} and
Lemma~\ref{lemgordon},
where we use $\Tr[Y_r^{2p}]\le r\llVert Y_r\rrVert ^{2p}$. This yields
\begin{eqnarray*}
\EE\llVert X\rrVert&\le& e^{1/2\alpha}\bigl\{2\sqrt{\bigl\lceil
\sigma^2\bigr\rceil+ \lceil\alpha\log n\rceil}+2\sqrt{2\lceil\alpha
\log n\rceil}\bigr\}
\\
&\le& e^{1/2\alpha}\{2\sigma+ 2\sqrt{\alpha\log n + 2} +2\sqrt{2\alpha
\log
n+2}\}
\end{eqnarray*}
for the choice $p = \lceil\alpha\log n\rceil$. If $n\ge2$ and
$\alpha\ge1$, then $2\le3\log2\le3\alpha\log n$, so
\[
\EE\llVert X\rrVert\le e^{1/2\alpha}\{2\sigma+ 6\sqrt{2\alpha\log n}\}.
\]
Defining $e^{1/2\alpha}=1+\varepsilon$ and noting that
$\varepsilon\le1/2$ implies $\alpha\ge1$ yields the result provided that
$n\ge2$ and $p\ge2$. The remaining cases are easily dealt with
separately. The result holds trivially in the case $n=1$. On the other
hand, the case $p=1$ can only occur when $\alpha\log n\le1$ and thus
$n\le2$. In this case we can estimate directly
$\EE\llVert X\rrVert \le\sqrt{n(\lceil\sigma^2\rceil+p)}\le
\sigma\sqrt{2}+2$ using Proposition~\ref{propmaincomp}.
\end{pf*}

%
\begin{rem}
Note that we use the moment method only to prove the comparison theorem
of Proposition~\ref{propmaincomp}; as will be seen below, this requires
only trivial combinatorics. All the usual combinatorial difficulties of
random matrix theory are circumvented by Lemma~\ref{lemgordon}, which
exploits the theory of Gaussian processes. After the first version of
this paper was posted, we learned from Shahar Mendelson that a related
idea has been used in \cite{Aub07} for a different purpose.
\end{rem}

%
\begin{rem}
The constant $6$ in the second term in Theorem~\ref{teomain} arises
from crude rounding in our proof. While this constant can be somewhat
improved for large $n$, our proof cannot yield a sharp constant here:
it can be verified in the example of the diagonal matrix
$b_{ij}=\mathbf{1}_{i=j}$ that the constant $\sqrt{2}$ in the
precise asymptotic $\EE\llVert X\rrVert \sim\sqrt{2\log n}$ cannot be
recovered
from our general proof. We therefore do not insist on optimizing this
constant, but rather state the convenient bound in
Theorem~\ref{teomain} which holds for any $n$. In contrast to the
constant in the second term, it was shown in the \hyperref[sec1]{Introduction} that the
constant in the first term is sharp.
\end{rem}

We now turn to the proof of Proposition~\ref{propmaincomp}. Let us begin
by recalling some standard observations. The quantity $\EE\Tr[X^{2p}]$
can be expanded as
\[
\EE\Tr\bigl[X^{2p}\bigr] = \sum_{u_1,\ldots,u_{2p}\in[n]}
b_{u_1u_2}b_{u_2u_3}\cdots b_{u_{2p}u_1} \EE[g_{u_1u_2}g_{u_2u_3}
\cdots g_{u_{2p}u_1}]. %
\]
Let $G_n=([n],E_n)$ be the complete graph on $n$ points, that is,
$E_n=\{\{u,u'\}\dvtx\break  u,u'\in[n]\}$. (Note that we have included self-loops.)
We will identify any $\mathbf{u}=(u_1,\ldots,u_{2p})\in[n]^{2p}$
with a
cycle $u_1\to u_2\to\cdots\to u_{2p}\to u_1$ in $G_n$ of length~$2p$.
If we denote by $n_i(\mathbf{u})$ the number of distinct edges that are
visited precisely $i$ times by the cycle $\mathbf{u}$, then we can write
[here $g\sim N(0,1)$]
\[
\EE\Tr\bigl[X^{2p}\bigr] = \sum_{\mathbf{u}\in[n]^{2p}}
b_{u_1u_2}b_{u_2u_3}\cdots b_{u_{2p}u_1} \prod
_{i\ge1}\EE\bigl[g^{i}\bigr]^{n_i(\mathbf{u})}.
\]
A cycle $\mathbf{u}$ is called \emph{even} if it visits each distinct edge
an even number of times, that is, if $n_i(\mathbf{u})=0$ whenever $i$ is
odd. As $\EE[g^i]=0$ when $i$ is odd, it follows immediately that the sum
in the above expression can be restricted to even cycles.

The \emph{shape} $\mathbf{s}(\mathbf{u})$ of a cycle $\mathbf{u}$ is
obtained by relabeling the vertices in order of appearance. For example,
the cycle $7\to3\to5\to4\to3\to5\to4\to3\to7$
has shape $1\to2\to3\to4\to2\to3\to4\to2\to1$.
We denote by
\[
\mathcal{S}_{2p}:=\bigl\{\mathbf{s}(\mathbf{u})\dvtx \mathbf{u}\mbox{
is an
even cycle of length }2p\bigr\} %
\]
the collection of shapes of even cycles, and we define the collection of
even cycles with given shape $\mathbf{s}$ starting (and ending) at a given
point $u$ as
\[
\Gamma_{\mathbf{s},u} := \bigl\{\mathbf{u}\in[n]^{2p}\dvtx \mathbf{s}(
\mathbf{u})=\mathbf{s}, u_1=u\bigr\} %
\]
for any $u\in[n]$ and $\mathbf{s}\in\mathcal{S}_{2p}$. Clearly the edge
counts $n_i(\mathbf{u})$ depend only on the shape $\mathbf{s}(\mathbf{u})$
of $\mathbf{u}$, and we can therefore unambiguously write
$n_i(\mathbf{s})$ for the number of distinct edges visited $i$ times by
any cycle with shape $\mathbf{s}$. We then obtain
\[
\EE\Tr\bigl[X^{2p}\bigr] = \sum_{u\in[n]} \sum
_{\mathbf{s}\in\mathcal{S}_{2p}} \prod_{i\ge1}\EE
\bigl[g^{i}\bigr]^{n_i(\mathbf{s})} \sum_{\mathbf{u}\in\Gamma_{\mathbf{s},u}}
b_{u_1u_2}b_{u_2u_3}\cdots b_{u_{2p}u_1}. %
\]
Finally, given any shape $\mathbf{s}=(s_1,\ldots,s_{2p})$, we denote by
$m(\mathbf{s})=\max_is_i$ the number of distinct vertices visited by any
cycle with shape $\mathbf{s}$.

Now that we have set up a convenient bookkeeping device, the proof of
Proposition~\ref{propmaincomp} is surprisingly straightforward.
It relies on two basic observations.

%
\begin{lem}
\label{lemcomp1}
Suppose that $\sigma_*\le1$. Then we have for any
$u\in[n]$ and $\mathbf{s}\in\mathcal{S}_{2p}$
\[
\sum_{\mathbf{u}\in\Gamma_{\mathbf{s},u}} b_{u_1u_2}b_{u_2u_3}\cdots
b_{u_{2p}u_1} \le\sigma^{2(m(\mathbf{s})-1)}. %
\]
In particular, it follows that
\[
\EE\Tr\bigl[X^{2p}\bigr] \le n \sum_{\mathbf{s}\in\mathcal{S}_{2p}}
\sigma^{2(m(\mathbf{s})-1)} \prod_{i\ge1}\EE
\bigl[g^{i}\bigr]^{n_i(\mathbf{s})}. %
\]
\end{lem}

\begin{pf}
Fix an initial point $u$, and shape $\mathbf{s}=(s_1,\ldots,s_{2p})$. Let
\[
i(k) = \inf\{j\dvtx s_j=k\} %
\]
for $1\le k\le m(\mathbf{s})$. That is, $i(k)$ is the first time in any
cycle of shape $\mathbf{s}$ at which its $k$th distinct vertex is visited
[of course, $i(1)=1$ by definition].

Now consider any cycle $\mathbf{u}\in\Gamma_{\mathbf{s},u}$.
As the cycle is even, the edge $\{u_{i(k)-1},u_{i(k)}\}$ must be visited
at least twice for every $2\le k\le m(\mathbf{s})$. On the other hand, as
the vertex $u_{i(k)}$ is visited for the first time at time $i(k)$, the
edge $\{u_{i(k)-1},u_{i(k)}\}$ must be distinct from the edges
$\{u_{i(\ell)-1},u_{i(\ell)}\}$ for all $\ell<k$. We can therefore
estimate
\begin{eqnarray*}
\sum_{\mathbf{u}\in\Gamma_{\mathbf{s},u}} b_{u_1u_2}b_{u_2u_3}\cdots
b_{u_{2p}u_1} &\le&\sum_{\mathbf{u}\in\Gamma_{\mathbf{s},u}} b_{uu_{i(2)}}^2
b_{u_{i(3)-1}u_{i(3)}}^2 \cdots b_{u_{i(m(\mathbf{s}))-1}u_{i(m(\mathbf{s}))}}^2
\\
&=& \sum_{v_2\ne\cdots\ne v_{m(\mathbf{s})}} b_{uv_2}^2
b_{v_{s_{i(3)-1}}v_3}^2 \cdots b_{v_{s_{i(m(\mathbf{s}))-1}}v_{m(\mathbf
{s})}}^2,
\end{eqnarray*}
where we use that $\max_{ij}\llvert b_{ij}\rrvert =\sigma_*\le1$.
As $s_{i(k)-1}<k$ by construction, it is readily seen that the quantity on
the right-hand side is bounded by $\sigma^{2(m(\mathbf{s})-1)}$.
\end{pf}

%
\begin{lem}
\label{lemcomp2}
Let $Y_r$ be defined as in Proposition~\ref{propmaincomp}.
Then for any $r>p$
\[
\EE\Tr\bigl[Y_r^{2p}\bigr] = r\sum
_{\mathbf{s}\in\mathcal{S}_{2p}} (r-1) (r-2)\cdots\bigl(r-m(\mathbf
{s})+1\bigr) \prod
_{i\ge1}\EE\bigl[g^{i}\bigr]^{n_i(\mathbf{s})}.
\]
\end{lem}

\begin{pf}
In complete analogy with the identity for $\EE\Tr[X^{2p}]$, we can write
\[
\EE\Tr\bigl[Y_r^{2p}\bigr] = \sum
_{\mathbf{s}\in\mathcal{S}_{2p}} \bigl\llvert\bigl\{\mathbf{u}\in
[r]^{2p}\dvtx
\mathbf{s}(\mathbf{u})=\mathbf{s}\bigr\}\bigr\rrvert\prod
_{i\ge1}\EE\bigl[g^{i}\bigr]^{n_i(\mathbf{s})}.
\]
Each cycle $\mathbf{u}\in[r]^{2p}$ with given shape
$\mathbf{s}(\mathbf{u})=\mathbf{s}$ is uniquely defined by
specifying its
$m(\mathbf{s})$ distinct vertices. Thus as long as $m(\mathbf{s})\le r$,
there are precisely
\[
r(r-1)\cdots\bigl(r-m(\mathbf{s})+1\bigr) %
\]
such cycles.
However, note that any even cycle of length $2p$ can visit at most
$m(\mathbf{s})\le p+1$ distinct vertices, so the assumption $p<r$ implies
the result.
\end{pf}

We can now complete the proof.

\begin{pf*}{Proof of Proposition~\ref{propmaincomp}}
Fix $p\in\mathbb{N}$, and let $r=\lceil\sigma^2\rceil+p$.
Then
\[
(r-1) (r-2)\cdots\bigl(r-m(\mathbf{s})+1\bigr) \ge\bigl(\sigma^2+p-m(
\mathbf{s})+1\bigr)^{m(\mathbf{s})-1} \ge\sigma^{2(m(\mathbf{s})-1)} %
\]
for any $\mathbf{s}\in\mathcal{S}_{2p}$, where we have used that
any even cycle of length $2p$ can visit at most
$m(\mathbf{s})\le p+1$ distinct vertices. It remains to apply
Lemmas~\ref{lemcomp1} and \ref{lemcomp2}.
\end{pf*}

\section{Extensions and complements}\label{secexcom}

\subsection{Nonsymmetric matrices}\label{secrect}

Let $X$ be the $n\times m$ random rectangular matrix with
$X_{ij}=g_{ij}b_{ij}$, where $\{g_{ij}\dvtx1\le i\le n,1\le j\le m\}$ are
independent $N(0,1)$ random variables and $\{b_{ij}\dvtx1\le i\le
n,1\le
j\le
m\}$ are given scalars. While this matrix is not symmetric, one can
immediately obtain a bound on $\EE\llVert X\rrVert $ from Theorem~\ref
{teomain} by
applying the latter to the symmetric matrix
\[
\tilde X = \lleft[ \matrix{0 & X
\cr
X^* & 0} \rright]. %
\]
Indeed, it is readily seen that $\llVert \tilde X\rrVert =\llVert
X\rrVert $, so we obtain
\[
\EE\llVert X\rrVert\le(1+\varepsilon)\biggl\{2(\sigma_1\vee
\sigma_2) + \frac{6}{\sqrt{\log(1+\varepsilon)}} \sigma_*\sqrt{\log
(n+m)}\biggr\}
\]
for any $0<\varepsilon\le1/2$ with
\[
\sigma_1 := \max_i\sqrt{\sum
_j b_{ij}^2},\qquad
\sigma_2 := \max_j\sqrt{\sum
_i b_{ij}^2},\qquad\sigma_* := \max
_{ij}\llvert b_{ij}\rrvert. %
\]
While this result is largely satisfactory, it does not lead to a sharp
constant in the first term: it is known from asymptotic
theory \cite{Gem80} that when $b_{ij}=1$ for all $i,j$, we
have $\EE\llVert X\rrVert \sim\sqrt{n}+\sqrt{m}$ as $n,m\to\infty$ with
$n/m\to\gamma\in\,]0,\infty[$, while the above bound can
only give
the weaker inequality $\EE\llVert X\rrVert \le2(1+o(1))(\sqrt{n}\vee
\sqrt{m})$.
The latter bound can therefore be off by as much as a
factor $2$.

We can regain the lost factor and also improve the logarithmic term by
exploiting explicitly the bipartite structure of $\tilde X$ in the proof
of Theorem~\ref{teomain}. This leads to the following sharp analogue of
Theorem~\ref{teomain} for rectangular random matrices.

%
\begin{teo}
\label{teorect}
Let $X$ be the $n\times m$ matrix with $X_{ij}=g_{ij}b_{ij}$. Then
\[
\EE\llVert X\rrVert\le(1+\varepsilon)\biggl\{\sigma_1+
\sigma_2 + \frac{5}{\sqrt{\log(1+\varepsilon)}} \sigma_*\sqrt{\log
(n\wedge m)}\biggr\}
\]
for any $0<\varepsilon\le1/2$.
\end{teo}

As the proof of this result closely follows the proof of
Theorem~\ref{teorect}, we will only sketch the necessary
modifications to the
proof in the rectangular setting.

\begin{pf*}{Sketch of proof}
Let $G_{n,m}=([n]\sqcup[m],E_{n,m})$ be the complete bipartite graph
whose left and right vertices are indexed by $[n]$ and $[m]$, respectively
(i.e., with edges $E_{n,m}=\{(u,v)\dvtx u\in[n],v\in[m]\}$). We begin by
noting that
\begin{eqnarray*}
&& \EE\Tr\bigl[\bigl(XX^*\bigr)^p\bigr]
\\
&&\qquad = \sum_{\mathbf{u}\in[n]^p} \sum_{\mathbf{v}\in[m]^p}
b_{u_1v_1}b_{u_2v_1}b_{u_2v_2}b_{u_3v_2}\cdots
b_{u_pv_p}b_{u_1v_p} \prod_{i\ge1}\EE
\bigl[g^i\bigr]^{n_i(\mathbf{u},\mathbf{v})},
\end{eqnarray*}
where we denote by $n_{i}(\mathbf{u},\mathbf{v})$ the number of distinct
edges in $G_{n,m}$ that are visited precisely $i$ times by the cycle
$u_1\to v_1\to u_2\to v_2\to\cdots\to u_p\to v_p\to u_1$.
In direct analogy with the symmetric case, we can define the collection
$\mathcal{S}_{2p}$ of shapes of even cycles of length $2p$, and by
$\Gamma_{\mathbf{s},u}$ the collection of cycles with given shape
$\mathbf{s}\in\mathcal{S}_{2p}$ starting at a given point $u\in[n]$.
We denote by $n_i(\mathbf{s})$ the number of distinct edges that are
visited precisely $i$ times by $\mathbf{s}$, and by $m_1(\mathbf{s})$ and
$m_2(\mathbf{s})$ the number of distinct right and left vertices,
respectively, that are visited by $\mathbf{s}$ (i.e., the number of
distinct vertices that appear in even and odd positions in the cycle).

It is now straightforward to adapt the proofs of Lemmas~\ref{lemcomp1}
and \ref{lemcomp2} to the present setting. Assuming that $\sigma
_*\le
1$, the analogue of Lemma~\ref{lemcomp1} yields
\[
\EE\Tr\bigl[\bigl(XX^*\bigr)^{p}\bigr] \le n \sum
_{\mathbf{s}\in\mathcal{S}_{2p}} \sigma_1^{2m_1(\mathbf{s})}
\sigma_2^{2(m_2(\mathbf{s})-1)} \prod_{i\ge1}\EE
\bigl[g^{i}\bigr]^{n_i(\mathbf{s})}.
\]
On the other hand, let $Y_{r,r'}$ be the $r\times r'$ matrix whose entries
are independent $N(0,1)$ random variables.
Then the analogue of Lemma~\ref{lemcomp2} yields
\begin{eqnarray*}
&&\EE\Tr\bigl[\bigl(Y_{r,r'}Y_{r,r'}^*\bigr)^{p}
\bigr]
\\
&&\qquad =
r\sum_{\mathbf{s}\in\mathcal{S}_{2p}} (r-1)\cdots\bigl(r-m_2(
\mathbf{s})+1\bigr) r'\bigl(r'-1\bigr)\cdots
\bigl(r'-m_1(\mathbf{s})+1\bigr)
\\
&&\quad\qquad{}\times \prod
_{i\ge1}\EE\bigl[g^{i}\bigr]^{n_i(\mathbf{s})}
\end{eqnarray*}
when $r>p/2$ and $r'>p/2$. Choosing $r=\lceil\sigma_2^2+p/2\rceil$
and $r'=\lceil\sigma_1^2+p/2\rceil$ yields
\[
\EE\Tr\bigl[\bigl(XX^*\bigr)^{p}\bigr] \le\frac{n}{r}\EE\Tr
\bigl[\bigl(Y_{r,r'}Y_{r,r'}^*\bigr)^{p}\bigr]
\]
in analogy with Proposition~\ref{propmaincomp}.
To complete the proof, we note that adapting the argument of
Lemma~\ref{lemgordon} to the rectangular case (cf. \cite{DS01}) yields
\[
\EE\bigl[\llVert Y_{r,r'}\rrVert^{2p}\bigr]^{1/2p}
\le\sqrt{r} + \sqrt{r'} + 2\sqrt{p}. %
\]
We can therefore estimate (assuming without loss of generality
that $\sigma_*=1$)
\[
\EE\llVert X\rrVert\le\EE\Tr\bigl[\bigl(XX^*\bigr)^{p}
\bigr]^{1/2p} \le n^{1/2p}\Bigl\{{\textstyle\sqrt{\bigl
\lceil\sigma_1^2+p/2\bigr\rceil} +\sqrt{\bigl
\lceil\sigma_2^2+p/2\bigr\rceil}}+2\sqrt{p}\Bigr\}.
\]
Choosing $p=\lceil\alpha\log n\rceil$ and proceeding as in the proof of
Theorem~\ref{teomain} yields the result with a dimensional factor of
$\sqrt{\log n}$ rather than $\sqrt{\log(n\wedge m)}$. However, as
$\llVert X\rrVert =\llVert X^*\rrVert $, the latter bound follows
by exchanging the roles of $n$
and $m$.
\end{pf*}

\subsection{Non-Gaussian variables}\label{secnong}

We have phrased our main results in terms of Gaussian random matrices for
concreteness. However, note that the core argument of the proof of
Theorem~\ref{teomain}, the comparison principle of
Proposition~\ref{propmaincomp}, did not depend at all on the Gaussian
nature of the
entries: it is only subsequently in Lemma~\ref{lemgordon} that we
exploited the theory of Gaussian processes. The same observation applies
to the proof of Theorem~\ref{teorect}. As a consequence, we can
develop various extensions of our main results to more general
distributions of the entries.

Let us begin by considering the case of sub-Gaussian random variables.

%
\begin{cor}
\label{corsubg}
Theorems~\ref{teomain} and \ref{teorect}
remain valid if the independent Gaussian random variables $g_{ij}$ are
replaced by independent symmetric random variables $\xi_{ij}$ such that
$\EE[\xi_{ij}^{2p}]\le\EE[g^{2p}]$ for every $p\in\mathbb{N}$ and $i,j$
[$g\sim N(0,1)$].
\end{cor}

\begin{pf}
As $\xi_{ij}$ are assumed to be symmetric, $\EE[\xi_{ij}^p]=0$ when
$p$ is
odd. It therefore follows readily by inspection of the proof that
Proposition~\ref{propmaincomp} (and its rectangular
counterpart) remains valid under the present assumptions.
\end{pf}

Corollary~\ref{corsubg} implies, for example, that the conclusions of
Theorems~\ref{teomain} and~\ref{teorect} hold verbatim when
$g_{ij}$ are
replaced by independent Rademacher variables $\varepsilon_{ij}$, that is,
$\mathbf{P}[\varepsilon_{ij}=\pm1]=1/2$; see Section~\ref{secrad} below
for more on such matrices. The moment assumption
$\EE[\xi_{ij}^{2p}]\le\EE[g^{2p}]$ is somewhat\vspace*{2pt} unwieldy, however.
We can
obtain a similar result under standard sub-Gaussian tail assumptions.

%
\begin{cor}
\label{corsubg2}
If the independent Gaussian variables $g_{ij}$ are replaced by independent
random variables $\xi_{ij}$ that are centered and sub-Gaussian in the sense
\[
\EE[\xi_{ij}]=0,\qquad\mathbf{P}\bigl[\llvert\xi_{ij}
\rrvert\ge t\bigr] \le Ce^{-t^2/2c}\qquad\mbox{for all }t\ge0 \mbox{ and
}i,j, %
\]
then Theorems~\ref{teomain} and \ref{teorect} remain
valid up to a universal constant that depends on $C$ and $c$ only. That
is, we
have $\EE\llVert X\rrVert \lesssim\sigma+ \sigma_*\sqrt{\log n}$ in
the case of
Theorem~\ref{teomain}, and
$\EE\llVert X\rrVert \lesssim\sigma_1+\sigma_2+\sigma_*\sqrt{\log
(n\wedge m)}$
in the case of Theorem~\ref{teorect}.
\end{cor}

\begin{pf}
Let $X'$ be an independent copy of $X$. As $\EE X'=0$, we obtain by
Jensen's inequality $\EE\llVert X\rrVert =\EE\llVert X-\EE
X'\rrVert \le\EE\llVert X-X'\rrVert $.
The entries of the matrix $X-X'$ are still sub-Gaussian (with the constants
$C,c$ increasing by at most a constant factor), but are now symmetric as
well. We can therefore assume without loss of generality that $\xi_{ij}$
are\vspace*{2pt} symmetric sub-Gaussian random variables. Using the integration formula
$\EE[\xi^{2p}] = \int_0^\infty\mathbf{P}[\llvert \xi\rrvert \ge t^{1/2p}]\,dt$,
it is
readily shown that the random variables $\xi_{ij}/K$ satisfy
$\EE[\xi_{ij}^{2p}]\le\EE[g^{2p}]$ for all $p\in\mathbb{N}$,
where $K$ is
a constant that depends on $C,c$ only. The result follows from
Corollary~\ref{corsubg}.
\end{pf}

%
\begin{rem}
The main difference between Corollaries~\ref{corsubg} and \ref{corsubg2}
is that the bound of Corollary~\ref{corsubg2} is multiplied by an
additional constant factor as compared to Corollary~\ref{corsubg}.
Thus the constant in front of the leading term in
Corollary~\ref{corsubg2} is no longer sharp.
This is of little consequence in many applications (particularly in
nonasymptotic problems), but implies that we no longer capture the exact
asymptotics of Wigner matrices. The latter can sometimes be recovered at
the expense of increasing the logarithmic term in the estimate; see
Corollary~\ref{corsymm} below.
\end{rem}

The sub-Gaussian assumption of Corollary~\ref{corsubg} requires that the
random variables $\xi_{ij}$ have at worst Gaussian tails. For random
variables with heavier tails, the conclusion of Theorem~\ref{teomain}
cannot hold as stated. Consider, for example, the diagonal case where
$b_{ij}=\mathbf{1}_{i=j}$, so that $\sigma=\sigma_*=1$; then
$\llVert X\rrVert =\max_{i\le n}\llvert \xi_{ii}\rrvert $, which
must grow faster in the
heavy-tailed setting than the $\sim\sqrt{\log n}$ bound that would be
obtained if the conclusion of Theorem~\ref{teomain} were valid.
It seems reasonable to expect that in the case of heavy-tailed entries,
the $\sqrt{\log n}$ rate must be changed to a quantity that controls the
maximum of the heavy-tailed random variables under consideration. This is,
roughly speaking, the content of the following result. (We will work in the
setting of Theorem~\ref{teomain} for simplicity, though an entirely
analogous result can be proved in the setting of Theorem~\ref{teorect}.)

%
\begin{cor}
\label{corheavy}
Let $X$ be the $n\times n$ symmetric matrix with
$X_{ij}=\xi_{ij}b_{ij}$, where $\{\xi_{ij}\dvtx i\ge j\}$ are independent
centered random variables and $\{b_{ij}\dvtx i\ge j\}$ are
given scalars. If $\EE[\xi_{ij}^{2p}]^{1/2p}\le Kp^{\beta/2}$ for some
$K,\beta>0$ and all $p,i,j$, then
\[
\EE\llVert X\rrVert\lesssim\sigma+ \sigma_*\log^{(\beta\vee1)/2}n. %
\]
The universal constant in the above inequality depends on $K,\beta$ only.
\end{cor}

Let us note that as $\EE[\xi^{2p}]^{1/2p}\lesssim\sqrt{p}$ for
sub-Gaussian random variables $\xi$, Corollary~\ref{corheavy} reduces
to Corollary~\ref{corsubg2} in the sub-Gaussian setting.\vspace*{1pt} If we consider
subexponential random variables, for example, then
$\EE[\xi^{2p}]^{1/2p}\lesssim p$, and thus $\sqrt{\log n}$ must be replaced
by $\log n$ in the conclusion of Theorem~\ref{teomain}. These scalings
are precisely as expected, as the maximum of $n$ independent random
variables $\xi_i$ with $\EE[\xi_i^{2p}]^{1/2p}\sim p^{\beta/2}$ is of
order $\log^{\beta/2}n$. Note, however, that the logarithmic factor only
changes when the tails of the entries are heavier than Gaussian.
The $\sqrt{\log n}$ factor cannot be reduced, in general, when the entries
have lighter tails than Gaussian, as the universality property of many
random matrix models leads to essentially Gaussian behavior; see
Remark~\ref{remlogremoval} for further discussion and examples.

\begin{pf*}{Proof of Corollary \ref{corheavy}}
Symmetrizing as in the proof of Corollary~\ref{corsubg2}, we can assume
without loss of generality that $\xi_{ij}$ are symmetric random variables.
We will also assume without loss of generality that $\beta\ge1$, as the
case $\beta<1$ is covered by Corollary~\ref{corsubg2}.

Let $g_{ij}$ and $\tilde g_{ij}$ be i.i.d. $N(0,1)$ random variables,
and define $\eta_{ij} = g_{ij}\llvert \tilde g_{ij}\rrvert ^{\beta-1}$.
Then $\eta_{ij}$ are symmetric random variables, and
by Stirling's formula
\[
\EE\bigl[\eta_{ij}^{2p}\bigr]^{1/2p} = \biggl[
\frac{2^{p\beta}}{\pi} \Gamma\biggl(p+\frac{1}{2}\biggr)\Gamma\biggl
(p(\beta-1)+
\frac{1}{2}\biggr) \biggr]^{1/2p} \gtrsim p^{\beta/2}.
\]
If we denote by $\tilde X$ the matrix with
entries $\tilde X_{ij}=\eta_{ij}b_{ij}$, then it follows readily
from the trace identities in the proof of Theorem~\ref{teomain}
that
\[
\EE\Tr\bigl[X^{2p}\bigr] \le C^{2p} \EE\Tr\bigl[\tilde
X^{2p}\bigr] \le C^{2p}n\EE\llVert\tilde X\rrVert
^{2p} %
\]
for all $p$ and a universal constant $C$ (depending on $K,\beta$).
We therefore have
\[
\EE\llVert X\rrVert\lesssim\EE\bigl[\llVert\tilde X\rrVert^{2\lceil
\log n\rceil}
\bigr]^{
1/2\lceil\log n\rceil}. %
\]
Applying Theorem~\ref{teomain} conditionally
on the variables $\tilde g=\{\tilde g_{ij}\}$ yields
\[
\EE\bigl[\llVert\tilde X\rrVert^{2\lceil\log n\rceil}\mid\tilde
g\bigr] \le\biggl[
C\max_i\sqrt{\sum_jb_{ij}^2
\llvert\tilde g_{ij}\rrvert^{2\beta-2}} + C\sigma_*\max
_{ij}\llvert\tilde g_{ij}\rrvert^{\beta-1}
\sqrt{\log n} \biggr]^{2\lceil\log n\rceil} %
\]
for another universal constant $C$. (While the statement of
Theorem~\ref{teomain} only gives a bound on $\EE\llVert X\rrVert $, an
inspection of the\vspace*{1pt}
proof shows that what is in fact being bounded is the quantity
$\EE[\llVert X\rrVert ^{2\lceil\alpha\log n\rceil}]^{1/2\lceil\alpha
\log
n\rceil}$ with $\alpha=1/2\log(1+\varepsilon)$.)

We must now estimate the expectation of the right-hand side of this
equation. Note that $\Vert\max_{i\le n}\vert Z_i\vert\Vert_p \le
\EE[\sum_{i=1}^n\vert Z_i\vert^p]^{1/p} \le n^{1/p}\max_{i\le n}\Vert
Z_i\Vert_p$
for any random variables $Z_1,\ldots,Z_n$.
Using $\llVert \tilde g_{ij}\rrVert _{p}\lesssim\sqrt{p}$, a simple
computation
yields
\[
\Bigl\llVert\max_{ij}\llvert\tilde g_{ij}\rrvert
^{\beta-1} \Bigr\rrVert_{2\lceil\log n\rceil} \lesssim\log^{(\beta-1)/2}n.
\]
Similarly, using the Rosenthal-type inequality of \cite{BBLM05},
Theorem~8, we obtain
\[
\biggl\llVert\max_i\sum_jb_{ij}^2
\llvert\tilde g_{ij}\rrvert^{2\beta-2} \biggr\rrVert
_{\lceil\log n\rceil} \lesssim\sigma^2 + \sigma_*^2 \Bigl
\llVert\max_{ij}\llvert\tilde g_{ij}\rrvert
^{\beta-1}\Bigr\rrVert_{2\lceil
\log n\rceil}^2\log n. %
\]
Substituting these estimates into the above expression completes the
proof.
\end{pf*}

As was discussed above, the drawback of Corollary~\ref{corsubg2} and
\ref{corheavy} is that an additional universal constant is introduced as
compared to Theorem~\ref{teomain}. In the following result, we have
retained the sharp constant at the expense of a suboptimal scaling of the
logarithmic term: for example, when applied to Gaussian entries,
Corollary~\ref{corsymm} yields $\log n$ instead of $\sqrt{\log n}$ in
Theorem~\ref{teomain}. Nonetheless, Corollary~\ref{corsymm} can be
useful in that it captures the sharp asymptotics of the edge of the
spectrum of Wigner-type matrices as long as $\sigma$ dominates the
logarithmic term. Moreover, when the random variables $\xi_{ij}$ are
uniformly bounded, Corollary~\ref{corsymm} is sharper than
Corollary~\ref{corsubg} in that the leading term depends on the
variance rather than the uniform size of the entries; this will be
exploited in Section~\ref{sectail} below.

%
\begin{cor}
\label{corsymm}
Let $X$ be the $n\times n$ symmetric random matrix with
$X_{ij}=\xi_{ij}b_{ij}$, where $\{\xi_{ij}\dvtx i\ge j\}$ are independent
symmetric random variables with unit variance and $\{b_{ij}\dvtx i\ge
j\}$ are
given scalars. Then we have for any $\alpha\ge3$
\[
\EE\llVert X\rrVert\le e^{2/\alpha}\Bigl\{2\sigma+ 14\alpha\max
_{ij}\llVert\xi_{ij}b_{ij}\rrVert
_{2\lceil\alpha\log n\rceil} \sqrt{\log n}\Bigr\}. %
\]
\end{cor}

\begin{pf}
Let $\varepsilon_{ij}$
be i.i.d. Rademacher random variables independent of $X$, and
denote by $\tilde X$ the matrix with entries $\tilde X_{ij}=
\varepsilon_{ij}X_{ij}$. As we assumed that $\xi_{ij}$ are
symmetric random variables, evidently $X$ and $\tilde X$ have the same
distribution. We now apply Corollary~\ref{corsubg} to
$\tilde X$ conditionally on the matrix $X$. This yields
\[
\EE\llVert X\rrVert\le(1+\delta)\EE\biggl[2\sqrt{\max
_i\sum_j X_{ij}^2}
+ \frac{6}{\sqrt{\log(1+\delta)}} \max_{ij}\llvert X_{ij}\rrvert
\sqrt{\log n} \biggr] %
\]
for any $0<\delta\le1/2$.
We can estimate
\[
\EE\max_{ij}\llvert X_{ij}\rrvert\le\biggl[\sum
_{ij}\EE\llvert X_{ij}\rrvert
^{2\lceil\alpha\log n\rceil} \biggr]^{1/2\lceil\alpha\log n\rceil} \le
e^{1/\alpha} \max
_{ij}\llVert\xi_{ij}b_{ij}\rrVert
_{2\lceil\alpha\log n\rceil}. %
\]
On the other hand, by the Rosenthal-type inequality of \cite{BBLM05},
Theorem~8, we have
\[
\biggl\llVert\sum_j X_{ij}^2
\biggr\rrVert_{\lceil\alpha\log n\rceil} \le(1+\delta)\sum_{j}b_{ij}^2
+ \frac{2}{\delta} \Bigl\llVert\max_{j}X_{ij}^2
\Bigr\rrVert_{\lceil\alpha\log n\rceil} \lceil\alpha\log n\rceil, %
\]
so that
\[
\EE\biggl[\max_i\sum_j
X_{ij}^2 \biggr] \le e^{1/\alpha}\biggl\{ (1+\delta)
\sigma^2 + \frac{2e^{1/\alpha}}{\delta} \max_{ij}\llVert
\xi_{ij}b_{ij}\rrVert_{2\lceil\alpha\log n\rceil}^2 \lceil
\alpha\log n\rceil\biggr\}. %
\]
Choosing $\alpha$ such that $e^{1/\alpha}=1+\delta$ and
using $\delta\ge\log(1+\delta)$, the result follows by combining
the above
estimates and straightforward manipulations.
\end{pf}

\subsection{Dimension-free bounds}\label{secdimfree}

A drawback of the results obtained so far is that they depend explicitly
on the dimension $n$ of the random matrix. This dependence is sharp in
many natural situations; see Section~\ref{seclower} below. On the other
hand, the results of Lata{\l}a \cite{Lat05} and of Riemer and Sch\"utt
\cite{RS13} have shown that it is possible to obtain dimension-free
estimates, where $n$ is replaced by an ``effective dimension'' that is
defined in terms of a norm of the matrix of coefficients of the form
\[
\bigl\llvert(b_{ij})\bigr\rrvert_p := \biggl[\sum
_{ij} \llvert b_{ij}\rrvert^p
\biggr]^{1/p}. %
\]
While the bounds of \cite{Lat05,RS13} yield
suboptimal results in many cases, a dimension-free formulation has at
least two advantages. First, a low-dimen\-sional matrix can be embedded in
a high-dimensional space without changing its norm: for example, if all
$b_{ij}=0$ except $b_{11}=1$, then $\EE\llVert X\rrVert \sim1$, but
Theorem~\ref{teomain} yields a bound of order $\sqrt{\log n}$. The
advantage of
a dimension-free bound is that it automatically adapts to high-dimensional
matrices that possess approximate low-dimensional structure. Second,
dimension-free results can be used to study infinite-dimensional matrices,
while Theorem~\ref{teomain} is not directly applicable in this setting.

The following result provides a dimension-free analogue of
Theorems~\ref{teomain} and~\ref{teorect}. To prove it, we apply the
stratification technique developed in \cite{RS13}.

%
\begin{cor}
\label{cordimfree}
Let $X$ be the $n\times n$ symmetric matrix with
$X_{ij}=g_{ij}b_{ij}$, where $\{g_{ij}\dvtx i\ge j\}$ are i.i.d. $\sim N(0,1)$
and $\{b_{ij}\dvtx i\ge j\}$ are given scalars. Then
\[
\EE\llVert X\rrVert\lesssim\sigma+ \sigma_* \sqrt{\log\frac{\llvert
(b_{ij})\rrvert _p}{\sigma_*}}
\]
for any $1\le p<2$, where the universal constant depends on $p$ only.
Similarly, if $X$ is the $n\times m$ random rectangular matrix
$X_{ij}=g_{ij}b_{ij}$, then for any $1\le p<2$
\[
\EE\llVert X\rrVert\lesssim\sigma_1 + \sigma_2 +
\sigma_* \sqrt{\log\frac{\llvert (b_{ij})\rrvert _p}{\sigma_*}}. %
\]
\end{cor}

\begin{pf}
We will prove the result in the symmetric case; the proof in the
rectangular case is identical. We also assume without loss of generality
that $\sigma_*=1$.

Define the matrices $X^{(k)}$ for $k\ge1$ as
$X^{(k)}_{ij}=g_{ij}b_{ij}\mathbf{1}_{2^{-k}<\llvert b_{ij}\rrvert \le
2^{-k+1}}$,
so that
\[
\EE\llVert X\rrVert= \EE\biggl\llVert\sum_{k\ge1}X^{(k)}
\biggr\rrVert\le\EE\biggl\llVert\sum_{k<k_0}X^{(k)}
\biggr\rrVert+ \sum_{k\ge k_0}\EE\bigl\llVert
X^{(k)}\bigr\rrVert%
\]
for a constant $k_0$ to be chosen appropriately in the sequel.

Denote by $c(k)$ the number of nonzero entries of $X^{(k)}$. Then
\[
c(k) := \bigl\llvert\bigl\{ij\dvtx2^{-k}<\llvert b_{ij}\rrvert
\le2^{-k+1}\bigr\}\bigr\rrvert\le2^{kp}\bigl\llvert
(b_{ij})\bigr\rrvert_p^p. %
\]
Therefore, the nonzero entries of $X^{(k)}$ must be
contained in a submatrix of size $c(k)\times c(k)$. Applying
Theorem~\ref{teomain} to this submatrix yields
\[
\EE\bigl\llVert X^{(k)}\bigr\rrVert\lesssim2^{-k+1}\bigl\{
\sqrt{c(k)}+\sqrt{\log c(k)}\bigr\} \lesssim2^{-k(1-p/2)}\bigl\llvert
(b_{ij})\bigr\rrvert_p^{p/2}. %
\]
As $p<2$, the right-hand side decays geometrically.
Let $k_0$ be the smallest integer $k$ such that
$2^{-k(1-p/2)}\llvert (b_{ij})\rrvert _p^{p/2}\le1$. Then we can estimate
\[
\sum_{k\ge k_0} \EE\bigl\llVert X^{(k)}\bigr
\rrVert\lesssim\sum_{k\ge k_0} 2^{-k(1-p/2)}\bigl
\llvert(b_{ij})\bigr\rrvert_p^{p/2} \le\sum
_{k\ge0} 2^{-k(1-p/2)} \lesssim\sigma, %
\]
where we use $\sigma\ge\sigma_*=1$.
On the other hand, the matrix $\sum_{k<k_0}X^{(k)}$ has at most
\[
\sum_{k<k_0}c(k) \le\sum
_{k<k_0}2^{kp}\bigl\llvert(b_{ij})\bigr
\rrvert_p^p \lesssim2^{k_0p}\bigl\llvert
(b_{ij})\bigr\rrvert_p^p \lesssim\bigl\llvert
(b_{ij})\bigr\rrvert_p^{p+p^2/(2-p)} %
\]
entries by the definition of $k_0$. Applying Theorem~\ref{teomain}
completes the proof.
\end{pf}

Note that the scaling in Corollary~\ref{cordimfree} improves as we
increase $p$. Unfortunately, the constant blows up as $p\to2$, so we
need $p<2$ to obtain a nontrivial result.

%
\begin{rem}
\label{remdimfree}
Up to universal constants, the result of Corollary~\ref{cordimfree} is
strictly better than that of Theorems~\ref{teomain} and \ref{teorect}
as $\llvert (b_{ij})\rrvert _p \le n^{2/p}\sigma_*$. It greatly
improves the bounds
of \cite{RS13}. The bound of \cite{Lat05} is of a somewhat
different nature: Lata{\l}a proves the inequality
$\EE\llVert X\rrVert \lesssim\sigma_1+\sigma_2+\llvert
(b_{ij})\rrvert _4$. The latter bound is
significantly worse than Corollary~\ref{cordimfree} in most cases, but
they are not strictly comparable.

Let us emphasize, however, that all the notions of effective dimension
used here or in \cite{Lat05,RS13} are essentially ad-hoc constructions.
As will be shown in Section~\ref{seclower} below, the bounds of
Theorems~\ref{teomain} and \ref{teorect} are tight in situations where
the
coefficients $b_{ij}$ exhibit a sufficient degree of homogeneity. The
improvement provided by the dimension-free bounds is therefore of interest
only in those cases where there is significant inhomogeneity in the
magnitude of the coefficients, that is, in the presence of many scales.
It is, however, unreasonable to expect that such inhomogeneity can be
captured in a sharp manner by a norm of the form $\llvert
(b_{ij})\rrvert _p$. This is
already illustrated by the simplest Gaussian examples: for example, if
$b_{ii}=1/\sqrt{\log i}$ and $b_{ij}=0$ for $i\ne j$, then a standard
Gaussian computation shows that $\EE\llVert X\rrVert \lesssim1$,
while all
dimension-free bounds we have discussed grow at least as $\sqrt{\log n}$.
\end{rem}

\subsection{Tail bounds}\label{sectail}

Given explicit bounds on the expectation $\EE\llVert X\rrVert $, we
can readily
obtain nonasymptotic tail inequalities for $\llVert X\rrVert $ by
applying standard
concentration techniques. In view of the significant utility of such tail
inequalities in applications, we record some useful results along
these lines here.

%
\begin{cor}
\label{corconcmain}
Under the assumptions of Theorem~\ref{teomain}, we have
\[
\mathbf{P} \biggl[ \llVert X\rrVert\ge(1+\varepsilon)\biggl\{2\sigma+
\frac{6}{\sqrt{\log(1+\varepsilon)}} \sigma_*\sqrt{\log n}\biggr\} + t
\biggr] \le e^{-t^2/4\sigma_*^2}
\]
for any $0<\varepsilon\le1/2$ and $t\ge0$. In particular, for every
$0<\varepsilon\le1/2$ there exists a universal constant
$c_\varepsilon$
such that for every $t\ge0$
\[
\mathbf{P}\bigl[ \llVert X\rrVert\ge(1+\varepsilon)2\sigma+ t \bigr]
\le
ne^{-t^2/c_\varepsilon\sigma_*^2}. %
\]
\end{cor}

\begin{pf}
As $\llVert X\rrVert =\sup_v\llvert \langle v,Xv\rangle\rrvert $
(the supremum is
over the unit ball) and
\[
\EE\bigl[\langle v,Xv\rangle^2\bigr] = \sum
_ib_{ii}^2v_i^4
+ 2\sum_{i\ne j}b_{ij}^2v_i^2v_j^2
\le2\sigma_*^2, %
\]
the first inequality follows from Gaussian concentration
\cite{BLM13}, Theorem 5.8, and Theorem~\ref{teomain}. For the second
inequality, note that we can estimate
\begin{eqnarray*}
\mathbf{P}\bigl[\llVert X\rrVert\ge(1+\varepsilon)2\sigma+
c_\varepsilon
\sigma_* t\bigr] &\le& \mathbf{P}\bigl[\llVert X\rrVert\ge(1+\varepsilon
)2\sigma+
c_\varepsilon'\sigma_*\sqrt{\log n}+\sigma_*t\bigr]
\\
&\le& e^{-t^2/4}
\end{eqnarray*}
for $t\ge2\sqrt{\log n}$ (with $c_\varepsilon,c_\varepsilon'$ chosen
in the obvious manner), while
\[
\mathbf{P}\bigl[\llVert X\rrVert\ge(1+\varepsilon)2\sigma+
c_\varepsilon
\sigma_* t\bigr] \le1 \le ne^{-t^2/4} %
\]
for $t\le2\sqrt{\log n}$. Combining these bounds completes the proof.
\end{pf}

Tail bounds on $\llVert X\rrVert $ have appeared widely in the recent
literature under the name ``matrix concentration inequalities''; see
\cite{Oli10,Tro12}. In the present setting, the
corresponding result of this kind implies that for all $t\ge0$
\[
\mathbf{P}\bigl[\llVert X\rrVert\ge t\bigr] \le ne^{-t^2/8\sigma^2}. %
\]
The second inequality of Corollary~\ref{corconcmain}
was stated for comparison with this matrix concentration bound.
Unlike the matrix concentration bound, Corollary~\ref{corconcmain} is
essentially optimal in that it captures not only the correct mean, but
also the correct tail behavior of $\llVert X\rrVert $ \cite{LT91},
Corollary 3.2
(Note that due to the factor $1+\varepsilon$ in the leading term, we
do not expect to see Tracy--Widom fluctuations at this scale.)

%
\begin{rem}
Integrating the tail bound obtained by the matrix concentration method
yields the estimate $\EE\llVert X\rrVert \lesssim\sigma\sqrt{\log
n}$. This method
therefore yields an alternative proof of the noncommutative Khintchine
bound that was discussed in the \hyperref[sec1]{Introduction}. Combining this bound with
concentration as in the proof of Corollary~\ref{corconcmain} already
yields a better tail bound than the one obtained directly from the matrix
concentration method. Nonetheless, it should be emphasized that the
suboptimality of the above bound on the expected norm stems from the
suboptimal tail behavior obtained by the matrix concentration method. Our
sharp tail bounds help clarify the source of this inefficiency: the
parameter $\sigma$ should only control the mean of $\llVert X\rrVert
$, while the
fluctuations are controlled entirely by $\sigma_*$.
\end{rem}

An entirely analogous result can be obtained in the rectangular setting of
Theorem~\ref{teorect}. As the proof is identical, we simply state the
result.

%
\begin{cor}
\label{corconcrect}
Under the assumptions of Theorem~\ref{teorect}, we have
\[
\mathbf{P} \biggl[ \llVert X\rrVert\ge(1+\varepsilon)\biggl\{
\sigma_1+\sigma_2 + \frac{5}{\sqrt{\log(1+\varepsilon)}} \sigma_*\sqrt{
\log(n\wedge m)}\biggr\} + t \biggr] \le e^{-t^2/2\sigma_*^2} %
\]
for any $0<\varepsilon\le1/2$ and $t\ge0$. In particular, for every
$0<\varepsilon\le1/2$ there exists a universal constant
$c_\varepsilon'$
such that for every $t\ge0$
\[
\mathbf{P}\bigl[ \llVert X\rrVert\ge(1+\varepsilon) (\sigma_1+
\sigma_2) + t \bigr] \le(n\wedge m)e^{-t^2/c_\varepsilon'\sigma_*^2}. %
\]
\end{cor}

The Gaussian concentration property used above is specific to Gaussian
variables. However, there are many other situations where strong
concentration results are available \cite{BLM13}, and where similar
results can be obtained. For example, if the Gaussian variables $g_{ij}$
are replaced by symmetric random variables $\xi_{ij}$ with
$\llVert \xi_{ij}\rrVert _\infty\le1$ (this captures in particular
the case of
Rademacher variables), Corollaries~\ref{corconcmain} and
\ref{corconcrect} remain valid with slightly larger universal constants
$c_\varepsilon,c_\varepsilon'$. This follows from the identical
proof, up
to the replacement of Gaussian concentration by a form of Talagrand's
concentration inequality \cite{BLM13}, Theorem 6.10.

In the case of bounded entries, however, a more interesting question is
whether it is possible to obtain tail bounds that capture the variance of
the entries rather than their uniform norm (which is often much bigger
than the variance), akin to the classical Bernstein inequality for sums of
independent random variables. We presently develop a very useful result
along these lines.

%
\begin{cor}
\label{corconcvar}
Let $X$ be an $n\times n$ symmetric matrix whose entries $X_{ij}$ are
independent symmetric random variables. Then there exists for any
$0<\varepsilon\le1/2$ a universal constant $\tilde c_\varepsilon$
such that for every $t\ge0$
\[
\mathbf{P}\bigl[\llVert X\rrVert\ge(1+\varepsilon)2\tilde\sigma+t\bigr
] \le
ne^{-t^2/\tilde c_\varepsilon\tilde\sigma_*^2}, %
\]
where we have defined
\[
\tilde\sigma:= \max_i \sqrt{\sum
_j \EE\bigl[X_{ij}^2\bigr]}, \qquad
\quad\tilde\sigma_* := \max_{ij}\llVert X_{ij}
\rrVert_{\infty}. %
\]
\end{cor}

\begin{pf}
Let $X_{ij}=\tilde X_{ij}\EE[X_{ij}^2]^{1/2}$, so that
$\tilde X_{ij}$ have unit variance. Then
\[
\EE\llVert X\rrVert\le(1+\varepsilon)2\tilde\sigma+ C_\varepsilon
\tilde
\sigma_*\sqrt{\log n} %
\]
for a suitable constant $C_\varepsilon$ by
Corollary~\ref{corsymm}. On the other hand,
a form of Talagrand's
concentration inequality \cite{BLM13}, Theorem 6.10, yields
\[
\mathbf{P}\bigl[\llVert X\rrVert\ge\EE\llVert X\rrVert+ t\bigr] \le
e^{-t^2/c\tilde\sigma_*^2} %
\]
for all $t\ge0$, where $c$ is a universal constant.
The proof is completed by combining these bounds as in the
proof of Corollary~\ref{corconcmain}.
\end{pf}

Corollary~\ref{corconcvar} should be compared with the matrix Bernstein
inequality in \cite{Tro12}, which reads as follows in our setting
(we omit the explicit constants):
\[
\mathbf{P}\bigl[\llVert X\rrVert\ge t\bigr] \le n e^{-t^2/c(\tilde
\sigma^2+\tilde\sigma_*t)}. %
\]
While this result looks quite different at first sight than
Corollary~\ref{corconcvar}, the latter yields strictly better tail behavior
up to universal constants: indeed, note that
\[
e^{-t^2/c^2\tilde\sigma_*^2} \le e^{1-2t/c\tilde\sigma_*} \le
3e^{-2t^2/c(\tilde\sigma^2+\tilde\sigma_*t)} %
\]
using $2x-1\le x^2$. The discrepancy between these results is readily
explained. In our sharp bounds, the variance term $\tilde\sigma$ only
appears in the mean of $\llVert X\rrVert $ and not in the
fluctuations: the latter
only depend on the uniform parameter $\tilde\sigma_*$ and do not capture
the variance. A tail bound in terms of $\tilde\sigma$ and
$\tilde\sigma_*$ should therefore indeed be of Hoeffding type, as in
Corollary~\ref{corconcvar}, rather than of Bernstein type as might be
expected from the ``matrix concentration'' approach. Using a Bernstein
form of Talagrand's concentration inequality \cite{Mas00}, Theorem 3,
in the proof of Corollary~\ref{corconcvar} does not lead to any further
improvement in the present setting.

%
\begin{rem}
If $X_{ij}$ in Corollary~\ref{corconcvar} are only
assumed to centered (rather than symmetric), we can symmetrize as in
the proof of Corollary~\ref{corsubg2} to obtain
\[
\mathbf{P}\bigl[\llVert X\rrVert\ge(1+\varepsilon)2\sqrt{2}\tilde\sigma
+t\bigr]
\le ne^{-t^2/\tilde c_\varepsilon\tilde\sigma_*^2}. %
\]
Unfortunately, this results in an additional factor $\sqrt{2}$ in
the leading term, which is suboptimal for Wigner matrices. We do not
know whether it is possible, in general, to improve the constant
when the entries are not symmetrically distributed.

Corollary~\ref{corconcvar} (and Corollary~\ref{corsymm} which is used
in its proof) also admit direct analogues in the setting of rectangular
matrices. As the proofs are essentially identical to the ones given
above, we leave such extensions to the reader.
\end{rem}

\subsection{Lower bounds}\label{seclower}

The main results of this paper provide upper bounds on $\EE\llVert
X\rrVert $.
However, a trivial lower bound already suffices to establish the sharpness
of our upper bounds in many cases of interest, at least for Gaussian
variables.

%
\begin{lem}
\label{lemlower}
In the setting of Theorem~\ref{teomain}, we have
\[
\EE\llVert X\rrVert\gtrsim\sigma+ \EE\max_{ij}\llvert
b_{ij}g_{ij}\rrvert. %
\]
Similarly, in the setting of Theorem~\ref{teorect}
\[
\EE\llVert X\rrVert\gtrsim\sigma_1 + \sigma_2 + \EE
\max_{ij}\llvert b_{ij}g_{ij}\rrvert.
\]
\end{lem}

\begin{pf}
Let us prove the second inequality; the first inequality follows in a
completely analogous manner. As $\llVert X\rrVert \ge\max
_{ij}\llvert X_{ij}\rrvert $, it is
trivial that
\[
\EE\llVert X\rrVert\ge\EE\max_{ij}\llvert X_{ij}
\rrvert= \EE\max_{ij}\llvert b_{ij}g_{ij}
\rrvert. %
\]
On the other hand, as $\llVert X\rrVert \ge\max_i\llVert
Xe_i\rrVert $ ($\{e_i\}$ is the
canonical basis in $\mathbb{R}^n$),
\[
\EE\llVert X\rrVert\ge\max_i\EE\llVert Xe_i
\rrVert\gtrsim\max_i\EE\bigl[\llVert Xe_i
\rrVert^2\bigr]^{1/2} = \sigma_2. %
\]
Here we use the estimate
\[
\EE\bigl[\llVert Xe_i\rrVert^2\bigr] = \bigl(\EE\llVert
Xe_i\rrVert\bigr)^2 + \Var\llVert Xe_i
\rrVert\lesssim\bigl(\EE\llVert Xe_i\rrVert\bigr)^2,
\]
where $\Var\llVert Xe_i\rrVert \le\max_jb_{ji}^2\lesssim
\max_j\EE[b_{ji}\llvert g_{ji}\rrvert ]^2
\le
(\EE\llVert Xe_i\rrVert )^2$ by the Gaussian
Poincar\'e inequality \cite{BLM13}, Theorem 3.20.
Analogously, we obtain
\[
\EE\llVert X\rrVert\ge\max_i\EE\bigl\llVert
X^*e_i\bigr\rrVert\gtrsim\sigma_1. %
\]
Averaging these three lower bounds yields the conclusion.
\end{pf}

This simple bound shows that our main results are sharp as long there are
enough large coefficients $b_{ij}$. This is the content of the following
easy bound.

%
\begin{cor}
\label{corsimplelower}
In the setting of Theorem~\ref{teomain}, suppose that
\[
\bigl\llvert\bigl\{ij\dvtx\llvert b_{ij}\rrvert\ge c\sigma_*\bigr\}
\bigr
\rrvert\ge n^{\alpha} %
\]
for some constants $c,\alpha>0$. Then
\[
\EE\llVert X\rrVert\asymp\sigma+ \sigma_*\sqrt{\log n}, %
\]
where the universal constant in the lower bound depends on
$c,\alpha$ only. The analogous result holds in the setting of
Theorem~\ref{teorect}.
\end{cor}

\begin{pf}
Denote by $I$ the set of indices in the statement of the corollary.
Then
\[
\EE\max_{ij}\llvert b_{ij}g_{ij}\rrvert
\ge\EE\max_{ij\in I}\llvert b_{ij}g_{ij}
\rrvert\ge c\sigma_*\EE\max_{ij\in I}\llvert g_{ij}
\rrvert\gtrsim\sigma_*\sqrt{\log\llvert I\rrvert} \gtrsim\sigma
_*\sqrt{\log n},
\]
where we use a standard lower bound on the maximum of independent
$N(0,1)$ random variables. The proof is completed by applying
Lemma~\ref{lemlower}.
\end{pf}

For example, it follows that our main results are sharp as soon as every
row of the matrix contains at least one large coefficient, that is, with
magnitude of the same order as $\sigma_*$. This is the case is many
natural examples of interest, and in these cases our results are optimal
(up to the values of universal constants).
Of course, it quite possible that our bound is sharp even when the
assumption of Corollary~\ref{corsimplelower} fails: for example, in view
of Lemma~\ref{lemlower}, our bound is sharp whenever $\sigma_*\sqrt
{\log
n}\lesssim\sigma$ regardless of any other feature of the problem.

Corollary~\ref{corsimplelower} suggests that our main results can
fail to
be sharp when the sizes of the coefficients $b_{ij}$ are very
heterogeneous. If the matrix contains a few large entries and many small
entries, one could still obtain good bounds by splitting the matrix into
two parts and applying Theorem~\ref{teomain} to each part; this is the
idea behind the dimension-free bounds of Corollary~\ref{cordimfree}.
However, when there are many different scales with few coefficients at
each scale, such an approach cannot be expected to yield sharp results in
general; see Remark~\ref{remdimfree} for a simple example.

%
\begin{rem}
\label{remschutt}
An intriguing observation that was made in \cite{RS13} is that the trivial
lower bound $\EE\llVert X\rrVert \ge\EE\max_i\llVert Xe_i\rrVert $
appears to be surprisingly
sharp: we do not know of any example where the corresponding upper bound
\[
\EE\llVert X\rrVert\stackrel{?} {\lesssim}\EE\max
_i\llVert Xe_i\rrVert%
\]
fails. If such an inequality were to hold, the conclusion of
Theorem~\ref{teomain} would follow easily using Gaussian
concentration and a
simple union bound.
In fact, if this were the case, we could obtain an improvement of
Theorem~\ref{teomain} in the following manner; cf. \cite{Tal14},
Proposition~2.4.16. Note that, by Gaussian concentration,
\[
\mathbf{P} \Bigl[\max_{i}\bigl\{\llVert Xe_i
\rrVert-\EE\llVert Xe_i\rrVert\bigr\}>t \Bigr] \le\sum
_{k} e^{-t^2/2\max_jb_{kj}^2} = \sum_{k}
k^{-t^2/2 \max_{ij}b_{ij}^2\log i}. %
\]
Integrating this bound therefore gives
\begin{eqnarray*}
\EE\llVert X\rrVert  &\stackrel{?} {\lesssim}& \EE\max
_i\llVert Xe_i\rrVert
\\
&\le&\max_i\EE\llVert Xe_i\rrVert+ \EE\max
_i\bigl\{\llVert Xe_i\rrVert-\EE\llVert
Xe_i\rrVert\bigr\}
\\
&\lesssim&\sigma+ \max_{ij}\llvert b_{ij}\rrvert
\sqrt{\log i}
\end{eqnarray*}
which would yield a strict improvement on Theorem~\ref{teomain}. [Note
that there is no loss of generality in sorting the rows of $(b_{ij})$ to
minimize the term $\max_{ij}\llvert b_{ij}\rrvert \sqrt{\log i}$.] We
do not know
of any mechanism, however, that would give rise to such inequalities, and
it is possible that the apparent sharpness of the quantity
$\EE\max_i\llVert Xe_i\rrVert $ is simply due to the fact that it is
of the same
order as the bound of Theorem~\ref{teomain} in most natural examples.
Regardless, it does not appear that our method of proof could be adapted
to give rise to inequalities of this form.
\end{rem}

%
\begin{rem}
The conclusion of Corollary~\ref{corsimplelower} relies heavily
on the Gaussian nature of the entries. When the distributions of
the entries are bounded, for example, it is possible that our bounds are
no longer sharp. This issue will be discussed further in
Section~\ref{secrad} below in the context of Rademacher matrices.
\end{rem}

\section{Examples}\label{secex}

\subsection{Sparse random matrices}\label{secsparse}

In the section, we consider the special case of Theorem~\ref{teomain}
where the coefficients $b_{ij}$ can take the values zero or one only. This
is in essence a sparse counterpart of Wigner matrices in which a subset of
the entries has been set to zero. This rather general model covers many
interesting random matrix ensembles, including the case of random band
matrices where $b_{ij}=\mathbf{1}_{\llvert i-j\rrvert \le k}$ that
has been of
significant recent interest \cite{Kho08,Sod10,BGP14}.

Let us fix a matrix $(b_{ij})$ of $\{0,1\}$-valued coefficients.
We immediately compute
\[
\sigma^2 = k,\qquad\sigma_* = 1, %
\]
where $k$ is the maximal number of nonzero entries in any row of the
matrix of coefficients $(b_{ij})$. If we interpret $(b_{ij})$ as the
adjacency matrix of a graph on $n$ points, then $k$ is simply the maximal
degree of this graph. The following conclusion follows effortlessly from
our main results.

%
\begin{cor}
\label{corsparse}
Let $X$ be the $n\times n$ symmetric random matrix with
$X_{ij}=g_{ij}b_{ij}$, where $\{g_{ij}\}$ are independent $N(0,1)$
variables and $b_{ij}\in\{0,1\}$. Let $k$ be the maximal number of
nonzero entries in a row of $(b_{ij})$. Then
\[
\EE\llVert X\rrVert\asymp\sqrt{k} + \sqrt{\log n}, %
\]
provided that every row of $(b_{ij})$ has at least one nonzero entry.
\end{cor}

\begin{pf}
This is immediate from Theorem~\ref{teomain} and Lemma~\ref{lemlower}.
\end{pf}

%
\begin{rem}
If a row of $(b_{ij})$ is zero, then the corresponding column is zero as
well by symmetry. We therefore lose nothing by removing this row and
column, and we can apply Corollary~\ref{corsparse} to
the resulting lower-dimensional matrix. The assumption that every row of
$(b_{ij})$ is nonzero is therefore completely innocuous.
\end{rem}

Our bound evidently captures precisely the correct order of magnitude of
the spectral norm of sparse random matrices. It is possible to obtain
much sharper conclusions, however, from our main results. To motivate
this, let us first quote a result that is stated in \cite{BGP14},
Theorem~2.3, under weaker assumptions. (For simplicity, we adopt in the
remainder of this section the setting and notation of
Corollary~\ref{corsparse}.)

%
\begin{teo}
\label{teobulk}
Suppose each row of $(b_{ij})$ has exactly $k$ nonzero entries.
Then the empirical spectral distribution of $X/\sqrt{k}$ converges
to the semicircle law
\[
\frac{1}{n}\sum_{i=1}^n
\delta_{\lambda_i(X/\sqrt{k})} \stackrel{k\to\infty} {\Longrightarrow}
\frac{1}{2\pi}
\sqrt{4-x^2} \mathbf{1}_{x\in[-2,2]} \,dx, %
\]
provided that $k=o(n)$.
[Here $\lambda_1(X)\ge\cdots\ge\lambda_n(X)$ are the eigenvalues of
$X$.]
\end{teo}

Theorem~\ref{teobulk} shows that the bulk of the spectrum of $X$ behaves
precisely like that of a Wigner matrix under minimal assumptions. As the
semicircle distribution has support $[-2,2]$, one might assume that edge
of the spectrum will converge to~$2$. That this is the case for Wigner
matrices is a textbook result \cite{AGZ10}. In the present case, however,
we obtain a phase transition phenomenon.

%
\begin{cor}
\label{corphtr}
Suppose that each row of $(b_{ij})$ has exactly $k$ nonzero entries.
Then the following phase transition occurs as
$n\to\infty$:
\begin{itemize}
\item If $k/\log n\to\infty$, then $\llVert X\rrVert /\sqrt{k}\to2$
in probability.
\item If $k/\log n\to0$, then $\llVert X\rrVert /\sqrt{k}\to\infty$
in probability.
\item If $k\sim\log n$, then $\{\llVert X\rrVert /\sqrt{k}\}$ is bounded
but may not converge to $2$.
\end{itemize}
\end{cor}

\begin{pf}
If $k/\log n\to0$, then $\llVert X\rrVert /\sqrt{k}\ge
\max_{ij}\llvert g_{ij}b_{ij}\rrvert /\sqrt{k}$. As each row has a
nonzero entry,
the maximum is taken over at least $n/2$ independent $N(0,1)$ random
variables which is of order $\sqrt{2\log(n/2)}$ as $n\to\infty$.
Thus $\llVert X\rrVert /\sqrt{k}$ diverges.

For $k/\log n\to\infty$, we note that Corollary~\ref{corconcmain} yields
\[
\mathbf{P}\bigl[ \llVert X\rrVert/\sqrt{k}\ge2 + \varepsilon\bigr]
\le
ne^{-C_\varepsilon k} = n^{1-C_\varepsilon k/\log n} %
\]
for a suitable constant $C_\varepsilon$. Thus
$\llVert X\rrVert /\sqrt{k}\le2+\varepsilon+ o(1)$ for any
$\varepsilon>0$.
On the other hand, Theorem~\ref{teobulk}
implies that $\llVert X\rrVert /\sqrt{k}\ge2-\varepsilon-o(1)$ for any
$\varepsilon>0$.

If\vspace*{1pt} $k=a\log n$, Corollary~\ref{corconcmain} similarly yields
that $\mathbf{P}[\llVert X\rrVert /\sqrt{k} > C]\to0$ for a
sufficiently large
constant $C$, so $\{\llVert X\rrVert /\sqrt{k}\}$ is bounded. However,
Lemma~\ref{lemlower} shows that $\EE\llVert X\rrVert /\sqrt{k} \gtrsim
a^{-1/2}>3$ when $a$ is sufficiently small.
\end{pf}

%
\begin{rem}
We have not investigated the precise behavior of $\llVert X\rrVert
/\sqrt{k}$ for the
boundary case $k=a\log n$. The proof of Corollary~\ref{corphtr} shows
that if $a$ is chosen sufficiently small, the
rescaled norm remains bounded but strictly separated from the bulk as
$n\to\infty$. We do not know whether this is the case for all $a$, or
whether the norm does in fact converge to $2$ when $a$ is sufficiently
large. A precise investigation of this question is beyond the scope of
this paper.
\end{rem}

In the special case of band matrices, Corollary~\ref{corphtr} was proved
by Sodin \cite{Sod10} following an earlier suboptimal result of Khorunzhiy
\cite{Kho08}. However, his combinatorial proof relies on the specific
positions of the nonzero entries of $(b_{ij})$. In a recent paper,
Benaych-Georges and P\'ech\'e \cite{BGP14} showed in the general setting
(i.e., without assuming specific positions of the entries) that
$\llVert X\rrVert /\sqrt{k}\to2$ when $k/\log^{9} n\to\infty$. To the
best of our
knowledge, however, the result of Corollary~\ref{corphtr} is new. While
this result is of independent interest, we particularly
emphasize how effortlessly a sharp conclusion could be derived from the
main results of this paper.

Beyond the Gaussian case, Corollaries~\ref{corheavy} and \ref
{corsymm} can
be used to obtain similar results in the presence of heavy-tailed entries.\vspace*{1pt}
Using Corollary~\ref{corheavy}, it can be shown that $\llVert
X\rrVert /\sqrt{k}$
remains bounded if and only if $\log^{\beta/2} n = O(k)$ when the entries
$\xi_{ij}$ have moments of order $\EE[\xi_{ij}^{2p}]^{1/2p}\sim
p^{\beta/2}$ with $\beta\ge1$. This establishes the appropriate phase
transition point in the heavy-tailed setting; however, we cannot conclude
convergence to the edge of the semicircle due to the additional universal
constant in Corollary~\ref{corheavy}. On the other hand, using
Corollary~\ref{corsymm} we can establish convergence to the edge of
the semicircle
under an assumption on the rate of growth of $k$ that is suboptimal by a
logarithmic factor; this is comparable to the results in \cite{BGP14},
though we obtain a somewhat better scaling. The details are omitted.

\subsection{Rademacher matrices}\label{secrad}

We have seen that our main results provide sharp bounds in many cases on
the norm of matrices with independent Gaussian entries. While our upper
bounds continue to hold for sub-Gaussian variables, this is not the case
for the lower bounds in Section~\ref{seclower}, and in this case we
cannot expect our results to be sharp at the same level of generality. As
a simple example, consider the case where $X$ is the diagonal matrix with
i.i.d. entries on the diagonal. If the entries are Gaussian, then
$\llVert X\rrVert \gtrsim\sqrt{\log n}$, so that Theorem~\ref
{teomain} is
sharp. If
the entries are bounded, however, then $\llVert X\rrVert \lesssim1$.
On the other
hand, the universality property of Wigner matrices shows that
Theorem~\ref{teomain} is sharp in this case even when adapted to
bounded random
variables (Corollary~\ref{corsubg2}).

In view of these observations, it is natural to ask whether it is possible
to obtain systematic improvement of our main results that captures the
size of the norm of random matrices with bounded entries. For
concreteness, let us consider the case of Rademacher matrices
$X_{ij}=\varepsilon_{ij}b_{ij}$, where $\{\varepsilon_{ij}\}$ are
independent Rademacher (symmetric Bernoulli) random variables. In this
setting, we can immediately obtain a trivial but useful improvement on
Corollary~\ref{corsubg2}. (In the rest of this section, we will consider
symmetric matrices and universal constants for simplicity; analogous
results for rectangular matrices or with explicit constants are easily
obtained.)

%
\begin{cor}
\label{corrad}
Let $X$ be the $n\times n$ symmetric random matrix with
$X_{ij}=\varepsilon_{ij}b_{ij}$, where $\{\varepsilon_{ij}\}$ are
independent Rademacher variables. Then
\[
\EE\llVert X\rrVert\lesssim(\sigma+\sigma_*\sqrt{\log n})\wedge\llVert
B\rrVert
, %
\]
where $B:=(\llvert b_{ij}\rrvert )$ is the matrix of absolute values
of the coefficients.
\end{cor}

\begin{pf}
In view of Corollary~\ref{corsubg}, it suffices to show that
$\EE\llVert X\rrVert \le\llVert B\rrVert $. Note, however, that
this inequality even holds
pointwise: indeed,
\[
\llVert X\rrVert= \sup_v \sum
_{ij} \varepsilon_{ij}b_{ij}v_iv_j
\le\sup_v \sum_{ij} \llvert
b_{ij}v_iv_j\rrvert= \llVert B\rrVert,
\]
where the supremum is taken over the unit ball in $\mathbb{R}^n$.
\end{pf}

Corollary~\ref{corrad} captures two reasons why a Rademacher matrix can
have small norm: either it behaves like a Gaussian matrix with small norm,
or its norm is uniformly bounded due to the boundedness of the matrix
entries. This idea mirrors the basic ingredients in the general theory of
Bernoulli processes \cite{Tal14}, Chapter~5. While simple,
Corollary~\ref{corrad} captures at least the Wigner and diagonal examples
considered above, albeit in a somewhat ad-hoc manner.
We will presently show that a less trivial result can be easily derived
from Corollary~\ref{corrad} as well.

The norm of Rademacher matrices was first investigated in a general
setting by Seginer \cite{Seg00}. Using a delicate combinatorial method,
he proves in this case that $\EE\llVert X\rrVert \lesssim\sigma\log
^{1/4}n$. The
assumption of Rademacher entries is essential: that such a bound cannot
hold in the Gaussian case is immediate from the diagonal matrix example.
Let us show that this result is an easy consequence of
Corollary~\ref{corrad}.

%
\begin{cor}
\label{corseginer}
Let $X$ be the $n\times n$ symmetric random matrix with
$X_{ij}=\varepsilon_{ij}b_{ij}$, where $\{\varepsilon_{ij}\}$ are
independent Rademacher variables. Then
\[
\EE\llVert X\rrVert\lesssim\sigma\log^{1/4}n. %
\]
\end{cor}

\begin{pf}
Fix $u>0$. Let us split the matrix into two parts $X=X^+ + X^-$, where
$X^+_{ij}=\varepsilon_{ij}b_{ij}\mathbf{1}_{\llvert b_{ij}\rrvert
>u}$ and
$X^-_{ij}=\varepsilon_{ij}b_{ij}\mathbf{1}_{\llvert b_{ij}\rrvert \le u}$.
For $X^-$, we can estimate
\[
\EE\bigl\llVert X^-\bigr\rrVert\lesssim\sigma+ u\sqrt{\log n}. %
\]
On the other hand, we estimate for $X^+$ by the
Gershgorin circle theorem
\[
\EE\bigl\llVert X^+\bigr\rrVert\le\bigl\llVert\bigl(\llvert b_{ij}
\rrvert\mathbf{1}_{\llvert b_{ij}\rrvert >u}\bigr)\bigr\rrVert\le
\max_i
\sum_j \llvert b_{ij}\rrvert
\mathbf{1}_{\llvert b_{ij}\rrvert >u} \le\frac{\sigma^2}{u}. %
\]
We therefore obtain for any $u>0$
\[
\EE\llVert X\rrVert\lesssim\sigma+ u\sqrt{\log n} + \frac{\sigma^2}{u}.
\]
The proof is completed by optimizing over $u>0$.
\end{pf}

Corollary~\ref{corseginer} not only recovers Seginer's result with a much
simpler proof, but also effectively explains why the mysterious term
$\log^{1/4}n$ arises. More generally, the method of proof suggests how
Corollary~\ref{corrad} can be used efficiently: we should attempt to
split the matrix $X$ into two parts, such that one part is small by
Theorem~\ref{teomain}, and the other part is small uniformly. This idea
also arises in a fundamental manner in the general theory of Bernoulli
processes \cite{Tal14}. Unfortunately, it is generally not clear for a
given matrix how to choose the best decomposition.

%
\begin{rem}
\label{remlogremoval}
In view of Corollary~\ref{corsparse}, one might hope that
Corollary~\ref{corrad} (or a suitable adaptation of this bound) could yield
sharp results in the general setting of sparse random matrices.
The situation for Rademacher matrices turns out to be more
delicate, however. To see this, let us consider two illuminating
examples. In the following, let $k=\lceil\sqrt{\log n}\rceil$, and
assume for simplicity that $n/k$ is integer.

First, consider the block-diagonal matrix $X$ of the form
\[
X = \lleft[ \matrix{ X_1 & & & &
\cr
& X_2 & & 0
\cr
& & \cdot&
\cr
& 0 & & \cdot
\cr
& & & & X_{n/k}} \rright],
\]
where each $X_i$ is a $k\times k$ symmetric matrix with
independent Rademacher entries. Such matrices are considered by
Seginer in \cite{Seg00}, who shows by an elementary argument that
$\EE\llVert X\rrVert \sim\sqrt{\log n}$. Thus Theorem~\ref{teomain} already
yields a
sharp result (and, in particular, the logarithmic term in
Theorem~\ref{teomain} cannot be eliminated).

On the other hand, Sodin \cite{Sod09} shows that if $X$ is the
Rademacher matrix where the coefficient matrix $B$ is chosen to be a
realization of the adjacency matrix of a random $k$-regular graph, then
$\EE\llVert X\rrVert \sim\sqrt{k}\le\log^{1/4}n$ with high
probability. Thus in
this case
$\EE\llVert X\rrVert \sim\sigma$, and it appears that the logarithmic
term in
Theorem~\ref{teomain} is missing (evidently none of our bounds are
sharp in this
case).

Note, however, that in both these examples the parameters
$\sigma,\sigma_*,\llVert B\rrVert $ are identical: we have $\sigma
=\sqrt{k}$,
$\sigma_*=1$, and $\llVert B\rrVert =k$ (by the Perron--Frobenius theorem).
In particular, there is no hope that the norm of sparse Rademacher
matrices can be controlled using only the degree of the graph: the
\emph{structure} of the graph must come into play. It is an interesting
open problem to understand precisely what aspect of this structure
controls the norm of sparse Rademacher matrices. This question is
closely connected to the study of random 2-lifts of graphs in
combinatorics \cite{BL06}.
\end{rem}


\section*{Acknowledgments}
In April 2014, one of us (ASB) attended the workshop ``Mathematical
Physics meets Sparse Recovery'' at Oberwolfach where the question of
removing the logarithmic term in the noncommutative Khintchine
inequality was raised by Joel Tropp in a more general context. While we
ultimately took a very different approach, we thank Joel for
motivating our interest in such problems which led us to develop the
ideas reported here. We are grateful to Sasha Sodin, Nikhil Srivastava
and Shahar Mendelson for several interesting discussions and for
providing us with a number of helpful references.
Afonso~S. Bandeira would also
like to thank Oberwolfach for its hospitality during the above-mentioned
workshop.


%

\printaddresses

\begin{thebibliography}{25}
\bibitem{AGZ10}
%
\begin{bbook}[mr]
\bauthor{\bsnm{Anderson},~\bfnm{Greg~W.}\binits{G.~W.}},
\bauthor{\bsnm{Guionnet},~\bfnm{Alice}\binits{A.}} \AND
\bauthor{\bsnm{Zeitouni},~\bfnm{Ofer}\binits{O.}}
(\byear{2010}).
\btitle{An Introduction to Random Matrices}.
\bseries{Cambridge Studies in Advanced Mathematics}
\bvolume{118}.
\bpublisher{Cambridge Univ. Press},
\blocation{Cambridge}.
\bid{mr={2760897}}
\end{bbook}
%
\bptok{imsref}%
\endbibitem

\bibitem{Aub07}
%
\begin{barticle}[mr]
\bauthor{\bsnm{Aubrun},~\bfnm{Guillaume}\binits{G.}}
(\byear{2007}).
\btitle{Sampling convex bodies: A random matrix approach}.
\bjournal{Proc. Amer. Math. Soc.}
\bvolume{135}
\bpages{1293--1303 (electronic)}.
\bid{doi={10.1090/S0002-9939-06-08615-1}, issn={0002-9939}, mr={2276637}}
\end{barticle}
%
\bptok{imsref}%
\endbibitem

\bibitem{BY88}
%
\begin{barticle}[mr]
\bauthor{\bsnm{Bai},~\bfnm{Z.~D.}\binits{Z.~D.}} \AND
\bauthor{\bsnm{Yin},~\bfnm{Y.~Q.}\binits{Y.~Q.}}
(\byear{1988}).
\btitle{Necessary and sufficient conditions for almost sure
convergence of the largest eigenvalue of a {W}igner matrix}.
\bjournal{Ann. Probab.}
\bvolume{16}
\bpages{1729--1741}.
\bid{issn={0091-1798}, mr={0958213}}
\end{barticle}
%
\bptok{imsref}%
\endbibitem

\bibitem{BGP14}
%
\begin{barticle}[mr]
\bauthor{\bsnm{Benaych-Georges},~\bfnm{Florent}\binits{F.}} \AND
\bauthor{\bsnm{P{\'e}ch{\'e}},~\bfnm{Sandrine}\binits{S.}}
(\byear{2014}).
\btitle{Largest eigenvalues and eigenvectors of band or sparse random
matrices}.
\bjournal{Electron. Commun. Probab.}
\bvolume{19}
\bpages{no. 4, 9}.
\bid{doi={10.1214/ECP.v19-3027}, issn={1083-589X}, mr={3164751}}
\end{barticle}
%
\bptok{imsref}%
\endbibitem

\bibitem{BL06}
%
\begin{barticle}[mr]
\bauthor{\bsnm{Bilu},~\bfnm{Yonatan}\binits{Y.}} \AND
\bauthor{\bsnm{Linial},~\bfnm{Nathan}\binits{N.}}
(\byear{2006}).
\btitle{Lifts, discrepancy and nearly optimal spectral gap}.
\bjournal{Combinatorica}
\bvolume{26}
\bpages{495--519}.
\bid{doi={10.1007/s00493-006-0029-7}, issn={0209-9683}, mr={2279667}}
\end{barticle}
%
\bptok{imsref}%
\endbibitem

\bibitem{BBLM05}
%
\begin{barticle}[mr]
\bauthor{\bsnm{Boucheron},~\bfnm{St{\'e}phane}\binits{S.}},
\bauthor{\bsnm{Bousquet},~\bfnm{Olivier}\binits{O.}},
\bauthor{\bsnm{Lugosi},~\bfnm{G{\'a}bor}\binits{G.}} \AND
\bauthor{\bsnm{Massart},~\bfnm{Pascal}\binits{P.}}
(\byear{2005}).
\btitle{Moment inequalities for functions of independent random variables}.
\bjournal{Ann. Probab.}
\bvolume{33}
\bpages{514--560}.
\bid{doi={10.1214/009117904000000856}, issn={0091-1798}, mr={2123200}}
\end{barticle}
%
\bptok{imsref}%
\endbibitem

\bibitem{BLM13}
%
\begin{bbook}[mr]
\bauthor{\bsnm{Boucheron},~\bfnm{St{\'e}phane}\binits{S.}},
\bauthor{\bsnm{Lugosi},~\bfnm{G{\'a}bor}\binits{G.}} \AND
\bauthor{\bsnm{Massart},~\bfnm{Pascal}\binits{P.}}
(\byear{2013}).
\btitle{Concentration Inequalities: A~Nonasymptotic Theory of Independence, With a Foreword by Michel Ledoux}.
\bpublisher{Oxford Univ. Press},
\blocation{Oxford}.
\bid{doi={10.1093/acprof:oso/9780199535255.001.0001}, mr={3185193}}
\end{bbook}
%
\bptok{imsref}%
\endbibitem

\bibitem{DS01}
%
\begin{bincollection}[mr]
\bauthor{\bsnm{Davidson},~\bfnm{Kenneth~R.}\binits{K.~R.}} \AND
\bauthor{\bsnm{Szarek},~\bfnm{Stanislaw~J.}\binits{S.~J.}}
(\byear{2001}).
\btitle{Local operator theory, random matrices and {B}anach spaces}.
In \bbooktitle{Handbook of the Geometry of {B}anach Spaces, {V}ol. {I}}
\bpages{317--366}.
\bpublisher{North-Holland},
\blocation{Amsterdam}.
\bid{doi={10.1016/S1874-5849(01)80010-3}, mr={1863696}}
\end{bincollection}
%
\bptok{imsref}%
\endbibitem

\bibitem{FK81}
%
\begin{barticle}[mr]
\bauthor{\bsnm{F{\"u}redi},~\bfnm{Z.}\binits{Z.}} \AND
\bauthor{\bsnm{Koml{\'o}s},~\bfnm{J.}\binits{J.}}
(\byear{1981}).
\btitle{The eigenvalues of random symmetric matrices}.
\bjournal{Combinatorica}
\bvolume{1}
\bpages{233--241}.
\bid{doi={10.1007/BF02579329}, issn={0209-9683}, mr={0637828}}
\end{barticle}
%
\bptok{imsref}%
\endbibitem

\bibitem{Gem80}
%
\begin{barticle}[mr]
\bauthor{\bsnm{Geman},~\bfnm{Stuart}\binits{S.}}
(\byear{1980}).
\btitle{A limit theorem for the norm of random matrices}.
\bjournal{Ann. Probab.}
\bvolume{8}
\bpages{252--261}.
\bid{issn={0091-1798}, mr={0566592}}
\end{barticle}
%
\bptok{imsref}%
\endbibitem

\bibitem{Kho08}
%
\begin{bincollection}[mr]
\bauthor{\bsnm{Khorunzhiy},~\bfnm{Oleksiy}\binits{O.}}
(\byear{2008}).
\btitle{Estimates for moments of random matrices with {G}aussian elements}.
In \bbooktitle{S\'eminaire de Probabilit\'es {XLI}}.
\bseries{Lecture Notes in Math.}
\bvolume{1934}
\bpages{51--92}.
\bpublisher{Springer},
\blocation{Berlin}.
\bid{doi={10.1007/978-3-540-77913-1_3}, mr={2483726}}
\end{bincollection}
%
\bptok{imsref}%
\endbibitem

\bibitem{Lat05}
%
\begin{barticle}[mr]
\bauthor{\bsnm{Lata{\l}a},~\bfnm{Rafa{\l}}\binits{R.}}
(\byear{2005}).
\btitle{Some estimates of norms of random matrices}.
\bjournal{Proc. Amer. Math. Soc.}
\bvolume{133}
\bpages{1273--1282 (electronic)}.
\bid{doi={10.1090/S0002-9939-04-07800-1}, issn={0002-9939}, mr={2111932}}
\end{barticle}
%
\bptok{imsref}%
\endbibitem

\bibitem{LT91}
%
\begin{bbook}[mr]
\bauthor{\bsnm{Ledoux},~\bfnm{Michel}\binits{M.}} \AND
\bauthor{\bsnm{Talagrand},~\bfnm{Michel}\binits{M.}}
(\byear{1991}).
\btitle{Probability in {B}anach Spaces: Isoperimetry and Processes}.
\bseries{Ergebnisse der Mathematik und Ihrer Grenzgebiete (3) [Results
in Mathematics and Related Areas (3)]}
\bvolume{23}.
\bpublisher{Springer},
\blocation{Berlin}.
\bid{doi={10.1007/978-3-642-20212-4}, mr={1102015}}
\end{bbook}
%
\bptok{imsref}%
\endbibitem

\bibitem{Mas00}
%
\begin{barticle}[mr]
\bauthor{\bsnm{Massart},~\bfnm{Pascal}\binits{P.}}
(\byear{2000}).
\btitle{About the constants in {T}alagrand's concentration
inequalities for empirical processes}.
\bjournal{Ann. Probab.}
\bvolume{28}
\bpages{863--884}.
\bid{doi={10.1214/aop/1019160263}, issn={0091-1798}, mr={1782276}}
\end{barticle}
%
\bptok{imsref}%
\endbibitem

\bibitem{Oli10}
%
\begin{barticle}[mr]
\bauthor{\bsnm{Oliveira},~\bfnm{Roberto~Imbuzeiro}\binits{R.~I.}}
(\byear{2010}).
\btitle{Sums of random {H}ermitian matrices and an inequality by {R}udelson}.
\bjournal{Electron. Commun. Probab.}
\bvolume{15}
\bpages{203--212}.
\bid{doi={10.1214/ECP.v15-1544}, issn={1083-589X}, mr={2653725}}
\end{barticle}
%
\bptok{imsref}%
\endbibitem

\bibitem{Pis03}
%
\begin{bbook}[mr]
\bauthor{\bsnm{Pisier},~\bfnm{Gilles}\binits{G.}}
(\byear{2003}).
\btitle{Introduction to Operator Space Theory}.
\bseries{London Mathematical Society Lecture Note Series}
\bvolume{294}.
\bpublisher{Cambridge Univ. Press},
\blocation{Cambridge}.
\bid{doi={10.1017/CBO9781107360235}, mr={2006539}}
\end{bbook}
%
\bptok{imsref}%
\endbibitem

\bibitem{RS13}
%
\begin{barticle}[mr]
\bauthor{\bsnm{Riemer},~\bfnm{Stiene}\binits{S.}} \AND
\bauthor{\bsnm{Sch{\"u}tt},~\bfnm{Carsten}\binits{C.}}
(\byear{2013}).
\btitle{On the expectation of the norm of random matrices with
non-identically distributed entries}.
\bjournal{Electron. J. Probab.}
\bvolume{18}
\bpages{no. 29, 13}.
\bid{issn={1083-6489}, mr={3035757}}
\end{barticle}
%
\bptok{imsref}%
\endbibitem

\bibitem{RV10}
%
\begin{binproceedings}[mr]
\bauthor{\bsnm{Rudelson},~\bfnm{Mark}\binits{M.}} \AND
\bauthor{\bsnm{Vershynin},~\bfnm{Roman}\binits{R.}}
(\byear{2010}).
\btitle{Non-asymptotic theory of random matrices: Extreme singular values}.
In \bbooktitle{Proceedings of the {I}nternational {C}ongress of
{M}athematicians. {V}olume {III}}
\bpages{1576--1602}.
\bpublisher{Hindustan Book Agency},
\blocation{New Delhi}.
\bid{mr={2827856}}
\end{binproceedings}
%
\bptok{imsref}%
\endbibitem

\bibitem{Seg00}
%
\begin{barticle}[mr]
\bauthor{\bsnm{Seginer},~\bfnm{Yoav}\binits{Y.}}
(\byear{2000}).
\btitle{The expected norm of random matrices}.
\bjournal{Combin. Probab. Comput.}
\bvolume{9}
\bpages{149--166}.
\bid{doi={10.1017/S096354830000420X}, issn={0963-5483}, mr={1762786}}
\end{barticle}
%
\bptok{imsref}%
\endbibitem

\bibitem{Sod09}
%
\begin{barticle}[mr]
\bauthor{\bsnm{Sodin},~\bfnm{Sasha}\binits{S.}}
(\byear{2009}).
\btitle{The {T}racy--{W}idom law for some sparse random matrices}.
\bjournal{J. Stat. Phys.}
\bvolume{136}
\bpages{834--841}.
\bid{doi={10.1007/s10955-009-9813-2}, issn={0022-4715}, mr={2545551}}
\end{barticle}
%
\bptok{imsref}%
\endbibitem

\bibitem{Sod10}
%
\begin{barticle}[mr]
\bauthor{\bsnm{Sodin},~\bfnm{Sasha}\binits{S.}}
(\byear{2010}).
\btitle{The spectral edge of some random band matrices}.
\bjournal{Ann. of Math. (2)}
\bvolume{172}
\bpages{2223--2251}.
\bid{doi={10.4007/annals.2010.172.2223}, issn={0003-486X}, mr={2726110}}
\end{barticle}
%
\bptok{imsref}%
\endbibitem

\bibitem{Tal14}
%
\begin{bbook}[mr]
\bauthor{\bsnm{Talagrand},~\bfnm{Michel}\binits{M.}}
(\byear{2014}).
\btitle{Upper and Lower Bounds for Stochastic Processes: Modern
Methods and Classical Problems}.
\bseries{Ergebnisse der Mathematik und Ihrer Grenzgebiete. 3. Folge.
A Series of Modern Surveys in Mathematics [Results in Mathematics and
Related Areas. 3rd Series. A Series of Modern Surveys in Mathematics]}
\bvolume{60}.
\bpublisher{Springer},
\blocation{Heidelberg}.
\bid{doi={10.1007/978-3-642-54075-2}, mr={3184689}}
\end{bbook}
%
\bptok{imsref}%
\endbibitem

\bibitem{Tao12}
%
\begin{bbook}[mr]
\bauthor{\bsnm{Tao},~\bfnm{Terence}\binits{T.}}
(\byear{2012}).
\btitle{Topics in Random Matrix Theory}.
\bseries{Graduate Studies in Mathematics}
\bvolume{132}.
\bpublisher{Amer. Math. Soc.},
\blocation{Providence, RI}.
\bid{mr={2906465}}
\end{bbook}
%
\bptok{imsref}%
\endbibitem

\bibitem{Tro12}
%
\begin{barticle}[mr]
\bauthor{\bsnm{Tropp},~\bfnm{Joel~A.}\binits{J.~A.}}
(\byear{2012}).
\btitle{User-friendly tail bounds for sums of random matrices}.
\bjournal{Found. Comput. Math.}
\bvolume{12}
\bpages{389--434}.
\bid{doi={10.1007/s10208-011-9099-z}, issn={1615-3375}, mr={2946459}}
\end{barticle}
%
\bptok{imsref}%
\endbibitem

\bibitem{Ver12}
%
\begin{bincollection}[mr]
\bauthor{\bsnm{Vershynin},~\bfnm{Roman}\binits{R.}}
(\byear{2012}).
\btitle{Introduction to the non-asymptotic analysis of random matrices}.
In \bbooktitle{Compressed Sensing}
\bpages{210--268}.
\bpublisher{Cambridge Univ. Press},
\blocation{Cambridge}.
\bid{mr={2963170}}
\end{bincollection}
%
\bptok{imsref}%
\endbibitem
\end{thebibliography}
\end{document}